\documentclass[11pt]{article}

\usepackage{latexsym}
\usepackage{amsfonts}
\usepackage{amsmath}
\usepackage{amssymb}

\newtheorem{theorem}{Theorem}
\newtheorem{definition}[theorem]{Definition}
\newtheorem{lemma}[theorem]{Lemma}
\newtheorem{corollary}[theorem]{Corollary}
\newtheorem{proposition}[theorem]{Proposition}
\newtheorem{remark}[theorem]{Remark}
\newtheorem{notation}[theorem]{Notation}
\numberwithin{theorem}{section}

\newcommand{\R}{I\!\!R}
\newcommand{\N}{I\!\!N}
\newcommand{\esssup}{\rm{esssup}}

\begin{document}

\title{Time consistent convex Feller processes and non linear second order partial differential equations}

\author{Jocelyne Bion-Nadal \\
CMAP (UMR CNRS 7641)  \\
Ecole Polytechnique
F-91128 Palaiseau cedex \\
jocelyne.bion-nadal@cmap.polytechnique.fr}

\date{July 6, 2012 \\}

\maketitle

\begin{abstract}
This paper shows how the theory of dynamic risk measures 
 provides viscosity solutions to a family of  second-order parabolic partial differential equations, even in the  degenerate  case. 

First, motivated by the martingale problem  approach of Stroock and Varadhan,  we make use of the time consistency characterization for dynamic risk measures,  to construct  time consistent convex Markov processes. This is done in a  general setting in which compacity arguments cannot be used, and for which there does not always exist an optimal control. Second, we prove that these processes lead to  viscosity solutions for  semi linear second-order partial differential equations with convex generator.

Finally we give  an application to mathematical finance.  
We show that our results allow for the construction of a No Arbitrage Pricing Procedure
in the context of stochastic volatility. Within this approach, convexity takes into account liquidy risk.\\

\end{abstract}

\section{Introduction} 
\label{sec:intro}
Under regularity assumptions the value function associated to an optimal problem for a non degenerate controlled diffusion is the classical solution of a second-order parabolic partial differential equation of  Hamilton-Jacobi-Bellman (HJB) type. In the degenerate case the HJB  equation does not have a classical solution. Lions introduced the notion of viscosity solution, and established  the link between viscosity solutions of second-order parabolic partial differential equations and  stochastic optimal control problems \cite{PL1,PL2}. 
 We refer to the paper of Crandall {\it et al.} 
 \cite{GIL} and to the book of  Fleming  and Soner \cite{FSo} for results and references.  Another approach making use of Backward Stochastic Differential Equations has been introduced by Pardoux and Peng \cite{PP2}. Recent papers in this direction by Briand and Hu \cite{BH}, Lio and Ley \cite{Da}, Delbaen {\it et al} \cite{DHR} and  Richou \cite{Ri} focuse on the case of convex generators or extend results to unbounded terminal conditions. Usually controlled diffusion is considered from the point of view of stochastic differential equations. An alternative approach has been introduced by Strook and Varadhan \cite{SV1,SV2}. It consists in the construction, given $s$ in $\R^+$ and $x$ in $\R^n$, of a probability measure $P_{s,x}$ on the set of continuous paths ${\cal C}([0,\infty[, \R^n)$ solution to a martingale problem associated to the diffusion.  Starting  from the stochastic differential equations approach with controls taking values in a compact subset of $\R^n$, El Karoui  {\it et al}  
\cite{EKNJ} use the  martingale problem approach to prove the existence of an optimal control.\\
On the other hand, the theory of dynamic risk measures on a filtered probability space has been developped in recent years. In the case of a Brownian filtration, dynamic risk measures coincide with $g$-expectations introduced by Peng \cite{SP}. An important property for dynamic risk measures is time consistency. The time consistency property for dynamic risk measures is the analogue of   the Dynamic Programming Principle.  For sublinear dynamic risk measures time consistency has been characterized by   Delbaen \cite{D}. For general convex dynamic risk measures two different  characterizations  of time consistency have been given. One by Cheridito  {\it et al} 
\cite{CDK}, the other by Bion-Nadal  \cite{BN03}. This last  characterization of time consistency  is very usefull in order to construct time consistent dynamic procedures. \\
The main goals of the present paper are twofold. First, motivated by the martingale problem  approach of  Stroock and Varadhan \cite{SV1,SV2},  we make use of the time consistency characterization proved in \cite{BN03},  to construct  time consistent convex Markov processes. Second, we prove that these processes lead to  viscosity solutions for  semi linear second order partial differential equations of the type 
\begin{equation}
\left \{\parbox{12cm}{
\begin{eqnarray} 
\partial_u v(u,x)+ {\cal L} v(u,x)+f(t, x,\sigma^*(u,x)Dv(u,x))&=&0 \nonumber \label{edp}\;\text{on}\; [0,t[ \times \R^n\\
v(t,x)& = & h(x)\nonumber 
\end{eqnarray}
}\right.
\end{equation}
where \begin{equation}
 {\cal L}v(t,x) = \frac{1}{2} Tr(\sigma\sigma^*(t,x)D^2(v)(t,x)) 
\label{lag}
\end{equation}
and $f: \R^+ \times \R^n \times \R^d \rightarrow \R$ is a Borel-measurable function such that  
$f(t,x,.)$ is a proper convex function. 
Here (and in all this paper) $\sigma$ is a given continuous bounded function with values in $M^n(\R)$ such that,  for all $(t,x)$, the matrix $\sigma(t,x)\sigma^*(t,x)$ is definite positive. These results do not require uniform ellipticity, that is they even apply  in the degenerate case. 
In addition to these goals, we give  an application to mathematical finance.\\ 

More precisely, the content of the paper is the following. For all given $r$ in $\R^+$ and $y$ in $\R^n$, we construct a whole time consistent dynamic procedure.  For this we need  a stable set of equivalent probability measures and a local penalty satisfying the cocycle condition \cite{BN03}. 
We start with a  family of probability measures on the set of continuous paths   which are solutions to a martingale problem starting from $y$ at time $r$. For every  probability measure in this family, the canonical process $(X_t)_{r \leq t}$ is a Feller process.
We then enlarge this set of probability measures in order to obtain a stable set. 
Given a  multivalued Borel mapping $\Lambda$ we  construct in Section \ref{secFellpro} such a stable set of probability 
measures $(\tilde {\cal Q}_{r,y})_S({\Lambda})$ indexed by  a set of progressively measurable processes $\mu(s,\omega)$ such that $\mu(s,\omega)$ belongs to $\Lambda (s,X_s(\omega))$. This set of probability measures  can be considered as a set of controls.\\  
The  construction of a time consistent dynamic procedure also  relies  on the construction of penalties associated to every probability measure in the stable set. For this we define penalties of the shape that we have considered in previous papers \cite{BN03} and \cite{BN04}. This construction  is generic in the sense that every normalized time consistent dynamic risk measure on a Brownian filtration is a limit of procedures for which the penalties are of this shape \cite{DPR}. The penalty associated to the probability measure $(P^{ \mu}_{r,y})$  is
$$
\alpha_{s,t}(P^{\mu}_{r,y})=E_{P^{\mu}_{r,y}}( \int\limits_{s}^{t} g(u,X_u(\omega),\mu(u,\omega)) d u|{\cal B}^r_{s})
$$
Given $(r,y)$ we construct this way (Section \ref{sectionTCC}) a time consistent convex dynamic procedure  defined for all ${\cal F}_t$ measurable  variable essentially bounded from below $Y$ by 
$$
\Pi^{r,y}_{s,t}(Y)=\esssup_{P^{\mu}_{r,y} \in (\tilde {\cal Q}_{r,y})_S({\Lambda})} (E_{P^{ \mu}_{r,y}}(Y|{\cal B}^r_{s})- \alpha_{s,t}(P^{\mu}_{r,y}))
$$
Notice that  for given $(r,y)$, in contrast with standard control theory, we do not construct a single value $v(r,y)$, but a whole time consistent dynamic procedure $(\Pi^{r,y}_{s,t})$. \\
 Remark that the canonical process $X_t$ is not a Feller process for every probability measure in the stable set $(\tilde {\cal Q}^{a}_{r,y})_S({\Lambda})$. However making use of the Feller property for the probability measures generating this stable set, we are able to prove that our time consistent procedures have the following Feller property (Section \ref{secFeller}):
  For all    Borelian  map  $h$ on $ \R^n$ bounded from below, there is a lower semi-continuous  function $\tilde h$ on   $[0,t[ \times \R^n$  such that 
$$ \forall s \in [0,t[,\; \forall x \in \R^n,\;\Pi^{s,x}_{s,t}(h(X_t))=\tilde h(s,x)
$$
$$\forall 0 \leq r \leq s <t,\;\;\Pi^{r,y}_{s,t}(h(X_t))=\tilde h (s,X_{s})
$$
 In case where $h$ is continuous bounded from below, the map $\tilde h $ extended by  $\tilde h (t,x)=h(x)$ is lower semi-continuous on $[0,t]\times\R^n$. \\
This result is proved under the `` weak'' hypothesis that the restriction of $g$ to $\Lambda$ is a Caratheodory function with polynomial growth. 
Let $f$ be the conjugate of $g$ on $\Lambda$, i.e.
$$f(u,x,z)=\sup_{\lambda \in \Lambda(u,x)}(z^*\lambda-g(u,x,\lambda))$$ 
 Assuming that $h$ is continuous, bounded from below and that the restriction of $g$ to   $\Lambda$ 
 is semi-continuous we prove (Section \ref{secsuper}) that the  lower semi-continuous function $\tilde h$ is a viscosity  supersolution of (1). This result is a consequence of the time consistency property which is the analogue of the dynamic programming principle. Notice that in the usual control theory the dynamic programming principle is proved making use of the existence of an optimal control. Here we make use of a sufficient condition for time consistency that we proved in  \cite{BN03}. This provides a wide class of procedures which are time consistent by construction,  and thus a wide class of equations of type (1) for which we get an explicit  viscosity supersolution.\\
Assuming furthermore a linear growth condition on $\Lambda$ and a continuity hypothesis on $f$, we prove that  the function $\tilde h$ is a viscosity solution of (1) if it is continuous (Section \ref{subsecU}).  A sufficient condition for the continuity of $\tilde h$ is that (1) satisfies a comparison principle.\\
It is important to notice that all the results of the present paper are in the setting of a control set which depends on $(t,X_t(\omega))$ and which is not bounded. The constructed stable set of probability measures is not necessarily compact for the weak topology. Therefore, compacity arguments cannot be used, and there does not always exist an optimal control.\\
Finally, in Section \ref{stocvol} 
we also provide an application to mathematical finance which illustrates the two main points of this paper. We show that the introduction of the multivalued Borel mapping allows for a construction of a No Arbitrage Pricing Procedure 
in the context of stochastic volatility. In this approach, convexity takes into account liquidy risk.
\section{Time consistent sublinear procedures in Markovian framework}
\label{secFellpro}
The aim of this section is to  use the  recent results on time consistent dynamic  procedures to construct sublinear  Feller processes.  Recall that in \cite{BN03,BN04}, we have already constructed examples of both sublinear and convex time consistent    dynamic procedures.  However the processes constructed in these papers do not have any Markov property. Here we will make use of probability measures solution to a martingale problem and of  the characterization of the time consistency property for dynamic risk measures  proved in \cite{BN03}, to provide families of  time consistent sublinear (convex in Section \ref{sectionTCC}) procedures with a Markov, even Feller, property.
\subsection{The martingale problem} 
\label{subsec1}
 Let $\Omega$ (resp. $\Omega^r$,  $\Omega^r_t$) be the set of continuous functions $\varphi$ defined on $\R^+$,  (resp. $[r,\infty[$ , $[r,t]$) with values in $\R^n$. 
  Recall that $\Omega$,  $\Omega^r$ and  $\Omega^r_t$ endowed with the topology of uniform convergence on compacts subsets of $\R^+$ are Polish spaces, that is  complete separable metrizable spaces.  Denote $(X_t)_{t \in \R^+}$  the coordinate process.
For given $r$, for  $r \leq s \leq t$, denote  ${\cal B}^s_{t}$ (resp.  ${\cal B}^s$) the  $\sigma$-algebra on $\Omega^r$ generated by the functions $X_u\;\; s \leq u \leq t$ (resp. the functions $X_u\;\; s \leq u $). 
The notation $Y^*$ denotes the transpose of  $Y$. Recall the following results from \cite{SV1, SV2}.  \\
 Assume that $a:\R^+ \times \R^n \rightarrow M_n(\R)$  is  continuous bounded strictly elliptic (i.e. positive definite at each point). It is proved in  \cite{SV1}, Theorem 4.2 and Theorem 5.6  that for 
 every $(s,x_0) \in \R^+ \times \R^n$, there is a unique probability measure $Q^{a}_{s,x_0}$ on $(\Omega^s,{\cal B}^s)$ such that 
\begin{equation}
X^s_{ \theta}(t)=exp \{\theta^*(X(t)-x_0) -\frac{1}{2} \int_s^t \theta^* a(u,X(u)) \theta du \} 
\label{eqmar0}
\end{equation}
is a $(Q^{a}_{s,x_0},({\cal B}^s_t)_{t \geq s})$ martingale for all  $\theta$ in $\R^n$, and such that $Q^a_{s,x_0}\{X_s=x_0\}=1$. 
Furthermore the process $(X_t)_{t \geq s}$ is then a strong Feller process for  $(Q^{a}_{s,x_0},({\cal B}^s_t)_{t \geq s})$.  Following  \cite{SV1}, one says that the probability measure $Q^a_{s,x_0}$ is the  unique solution to the martingale problem  (\ref{eqmar0}) starting from $x_0$ at time $s$.\\
 Let $\sigma$ be a continuous bounded map $\sigma:\R^+ \times \R^n \rightarrow M_n(\R)$ such that for all $(t,x)$, $a(t,x)=\sigma(t,x)\sigma^*(t,x)$.  Let $\lambda$ be a borelian bounded map  $\lambda:\R^+ \times \R^n \rightarrow \R^n$.
 From Theorem 6.2 of  \cite{SV1} there is a unique probability measure $Q^{\sigma, \lambda}_{s,x_0}$ on $(\Omega^s,{\cal B}^s)$ solution to the martingale problem:
\begin{equation}
exp \{(\theta)^*(X(t)-x_0)-\int_s^t \theta^*(\sigma \lambda)(u,X(u))  du -\frac{1}{2} \int_s^t \theta^* a(u,X(u)) \theta du \} 
\label{eqmar}
\end{equation}
starting from $x_0$ at time $s$, i.e. such that $Q^{\sigma, \lambda}_{s,x_0}(\{X_s=x_0\})=1$. 
For all $s \leq t$,
\begin{equation}
(\frac{dQ^{\sigma, \lambda}_{s,x_0}}{dQ^{a}_{s,x_0}})_{B^s_t}= exp(\int_s^t (\lambda^*\sigma^{-1})) (u,X_u)dX_u-\frac{1}{2} \int_s^t(\lambda^* \lambda )(u,X_u)du)
\label{eqDD}
\end{equation}
the stochastic integral being computed with respect to $Q^{a}_{s,x_0}$.
  From  \cite{SV2} Theorem 7.1, the process $(X_t)_{t \geq s}$ is then a strong Feller process  for   $(Q^{a}_{s,x_0},({\cal B}^s_t)_{t \geq s})$.\\
Notice that assuming $a$ continuous bounded strictly elliptic the non negative square root of  $a$, $\sigma_a=a^{\frac{1}{2}}$ is always continuous bounded and satisfies $\sigma_a\sigma_a^*=a$. However we do not want to restrict to this specific choice for $\sigma$. In particular the financial application (Section \ref{stocvol}) corresponds to a matrix $\sigma$ which is not symetric.
\begin{notation}
Given $\omega$ in $\Omega^r$, given $0 \leq r <s$ denote $\pi^{r,s}$
the continuous projection from $\Omega^r$ onto  $\Omega^s$:
\begin{eqnarray}
\pi^{r,s}(\omega')(u)=\omega'(u),\;\forall s \leq u 
\label{eqi}
\end{eqnarray}
For all probability measure $Q$ on $(\Omega^r, {\cal B}^r)$, $Q  (\pi^{r,s})^{-1}$ denotes the probability measure on $(\Omega^s,{\cal B}^s)$ image of $Q$ by $\pi^{r,s}$.
\end{notation}
Recall from \cite{SV1}, Theorem 2.1, that given $r<s$,  every  probability measure on $(\Omega,{\cal B}^r)$ admits  a regular conditional probability 
distribution given ${\cal B}^r_s$.
\begin{lemma}
Let $\sigma$ be continuous bounded and $a=\sigma \sigma^*$ strictly elliptic. Let  $\lambda$ be borelian bounded. Let $(Q^{\sigma, \lambda}_{r,x_0})$ be the unique probability measure solution to  the martingale  problem starting from $x_0$ at time $r$  (equation (\ref{eqmar})). Given $s>r$, let $(Q^{\sigma, \lambda}_{r,x_0})_{s, \omega}$ be a regular conditional probability  of  $Q^{\sigma, \lambda}_{r,x_0}$ given $B^r_{s}$.  There is a  $Q^{\sigma, \lambda}_{r,x_0}$-null set N such that 
 $(Q^{\sigma, \lambda}_{r,x_0})_{s, \omega}(\pi^{r,s})^{-1} =Q^{\sigma, \lambda}_{s, X_s(\omega)}$ for all $\omega$ in $N^c$.\\
In particular for all Borelian map $f$ on $\R^k$ and  all $s \leq t_1\leq ... \leq t_k$,
\begin{equation}
(Q^{\sigma, \lambda}_{r,x_0})_{s, \omega}(f(X_{t_1},..., X_{t_k})=Q^{\sigma, \lambda}_{s, X_s(\omega)}(f(X_{t_1},..., X_{t_k}),\;\;\forall \omega \in N^c
\label{eqcp1}
\end{equation}
\label{lemmapc1}
\end{lemma}
It follows from the  definition of a regular conditional probability  that \\$(Q^{\sigma, \lambda}_{r,x_0})_{s, \omega} (\pi^{r,s})^{-1})(\{\omega'\;|X_s(\omega)=X_s(\omega')\})=1$. Notice also that from the definition of $\pi^{r,s}$ for all   Borelian function $f$ on $\R^k$ and all $s \leq t_1\leq ... \leq t_k$,
$(Q^{\sigma, \lambda}_{r,x_0})_{s, \omega} (\pi^{r,s})^{-1})(f(X_{t_1},..., X_{t_k})=(Q^{\sigma, \lambda}_{r,x_0})_{s, \omega}(f(X_{t_1},..., X_{t_k})$.
 It follows then from Theorem 3.1 of \cite{SV1}, that  there
 is a set $N$ such that $Q^{\sigma, \lambda}_{r,x_0}(N)=0$ and such that for all $\omega \in N^c$,  
  $(Q^{\sigma, \lambda}_{r,x_0})_{s, \omega} (\pi^{r,s})^{-1}$  solves the martingale problem (\ref{eqmar}) starting  from $X_s(\omega)$ at time $s$. By unicity of the solution to this martingale problem, it follows that 
\begin{equation}
\forall \omega \in  N^c, \;\;  (Q^{\sigma, \lambda}_{r,x_0})_{s, \omega} (\pi^{r,s})^{-1}= Q^{\sigma, \lambda}_{s, X_s(\omega)}
\label{eqprobcon}
 \end{equation}
$\square $ 

\subsection{Estimates of the moments of $X_t$}
\label{secestimates}

\begin{proposition}
For all $q\geq 1$, $A$, $B$,  and $t>0$, there is  $K >0$ such that for all $y$ such that $||y|| \leq C$ and $\sigma$ continuous bounded such that
\begin{equation}
\forall \theta \in \R^n,\;0< \theta^* a \theta \leq A ||\theta||^2,\; with \;a=\sigma \sigma^*
\label{eqstrellip0}
\end{equation}
for $||\sigma \lambda|| \leq B$, 
for all $0 \leq r \leq s\leq t$, 
\begin{equation}
E_{Q^{\sigma, \lambda}_{r,y}}(\sup_{s \leq u \leq t }(||X_t-X_{u}||^{2q}) \leq K (t-s)^q(||y||^{2q}+1)
\label{eqKry2}
\end{equation}
In the following $||a||\leq A$ means that $a$ satisfies (\ref{eqstrellip0}).
\label{propKryl}
\end{proposition}
{\bf Proof}
\begin{itemize}
\item Assume that $a$ is uniformly elliptic and satisfies (\ref{eqstrellip0}). Denote $a^{\frac{1}{2}}$ the non negative square root of $a$. From Corollary 3.2 of \cite{SV1}, there is a Brownian motion $\beta$ with respect to $Q^{\sigma, \lambda}_{r,y}$ such that 
\begin{equation}
X_t=y + \int_{r}^{t}a^{\frac{1}{2}}(u,X_u)d\beta(u)+\int_{r}^{t}\sigma \lambda(u,X_u)du
\label{eqKry1}
\end{equation}
From \cite{Kr}, II 5 Corollary 10, there is a   constant  $K>0$  depending only on $q$ $A$ $B$ and $t$ such that  equation (\ref{eqKry2}) 
is satisfied.
\item For general $a$ satisfying (\ref{eqstrellip0}),  Let $b=\sigma\lambda$.  Notice that given $r$ and $y$, $Q^{\sigma,\lambda}_{r,y}$ depends only on $a=\sigma \sigma^*$ and $b=\sigma \lambda$. Let $R^{a,b}=Q^{\sigma,\lambda}_{r,x}$.  Consider $a_j$ continuous uniformly elliptic such that the sequence  $a_j$ is uniformly bounded by $A$ and uniformly convergent to $a$ on compact spaces. From Theorem 9.2 of \cite{SV2} applied for given $b$, $r$ and $y$, $R^{a_j,b}$ converges weakly to $R^{a, b}$  as $j \rightarrow \infty$.  Given $r \leq s\leq t$, consider a sequence ${\cal T}_n$ of subdivisions of $[s,t]$ whose step tends to $0$ as $n \rightarrow \infty$. For all $m \in \N^*$, $\inf(\sup_{u \in {\cal T}_n}(||X_t-X_{u}||^{2q},m))$ is a continuous bounded function. Passing to the limit as $j \rightarrow \infty$, we obtain that for all $m>0$,
$$E_{Q^{\sigma,\lambda}_{r,y}}(\inf(\sup_{u \in {\cal T}_n}(||X_t-X_{u}||^{2q},m)) \leq K (t-s)^q(||y||^{2q}+1)$$
Using the monotone convergence theorem, we obtain that  (\ref{eqKry2}) is satisfied  for all $a$ continuous bounded satisfying (\ref{eqstrellip0}) and  $||\sigma \lambda|| \leq B$. 
\end{itemize} \hfill $\square$
\begin{corollary}
Let $m \geq 1$. Let $\epsilon>0$. There is $k_0$ depending only on $m$ $A$,$B$,$C$ and $t$ such that for all $k>k_0$, for all $\sigma$ continuous bounded such that $\sigma\sigma^*$ satisfies (\ref{eqstrellip0}), for $||\sigma\lambda|| \leq B$,  $||y|| \leq C$,and $r\leq s \leq t$, 
\begin{equation} 
E_{Q^{\sigma, \lambda}_{r,y}}(\sup_{r \leq u \leq t }(||X_{u}-y||^{m}1_{\{\sup_{r \leq u \leq t}||X_u-y|| \geq k\}})<\epsilon
\label{eqcontf}
\end{equation}
\label{corKry}
\end{corollary}
{\bf Proof} 
From Proposition \ref{propKryl}, for all $k > 0$,
\begin{eqnarray}
k E_{Q^{\sigma, \lambda}_{r,y}}(\sup_{r \leq u \leq t }(||X_{u}-y||^{m}1_{\{\sup_{r \leq u \leq t}||X_u-y|| \geq k\}}) \leq \nonumber\\
\leq  E_{Q^{\sigma, \lambda}_{r,y}}(\sup_{r \leq u \leq t }(||X_{u}-y||^{m+1}) \leq K(t-r)^{\frac{m+1}{2}}(C^{m+1}+1)
\end{eqnarray}   \hfill $\square$

\begin{proposition} 
\begin{enumerate}
\item
Let $A, K,B>0$.  Given $r$ in $\R^+$,  the set
 of probability measures  $\{Q^{\sigma, \lambda}_{r,y},\;\;||a|| \leq A,\; ||y|| \leq   K, \||\sigma \lambda|| \leq B\}$ on $(\Omega^r,{\cal B}^r)$ is  weakly relatively compact. 
\item
For given $r$ in $\R^+$, the map $x \rightarrow Q^{\sigma,\lambda}_{r,x}$ is continuous for the weak topology. 
\end{enumerate}
\label{lemmawrc}
\end{proposition}
{\bf Proof}
\begin{enumerate}
\item
It follows from equation (\ref{eqKry2}) that for $q \geq 1$, there is a constant $K_1>0$ such that for $||a|| \leq A,\; ||y|| \leq   K, \||\sigma \lambda|| \leq B$,
\begin{equation}
E_{Q^{\sigma, \lambda}_{r,y}}(\sup_{s \leq u \leq t }(||X_t-X_{u}||^{2q}) \leq K_1(t-s)^q
\label{eqKry1-1}
\end{equation}
From theorem 2.3 of  \cite{SV1},  the set of probability measures  $\{Q^{\sigma, \lambda}_{r,y},\;\;||a||\leq A,\; ||y|| \leq   K,  ||\sigma \lambda|| \leq B,\}$ is thus weakly relatively compact.
\item 
Let $x_n$ be a sequence with limit $x$.  The set  $\{Q^{\sigma, \lambda}_{r,x_n}, \; n \in \N\}$ is  weakly relatively compact. Any limit  point of the family  $Q^{\sigma, \lambda}_{r,x_n}$ solves the martingale problem for $(\sigma,\lambda)$ starting from $(r,x)$. By unicity of the solution of the martingale problem it follows that $Q^{\sigma, \lambda}_{r,x_n} \rightarrow  Q^{\sigma \lambda}_{r,x}$ weakly.
\end{enumerate}
$\square $


\subsection{Feller processes with continuous paths} 
\label{subsec1-3}
\begin{proposition}
Let $\sigma$ be continuous bounded and $a=\sigma \sigma^*$ strictly elliptic. Let  $\lambda$ be borelian bounded.  Let $t>0$. For given $\phi$ Borelian bounded, there is a   continuous bounded function $U^{\sigma, \lambda}_t(\phi):(s,x)\in [0,t[\times \R^n \rightarrow  U^{\sigma, \lambda}_t(\phi)(s,x)$ such that:\\
For all  $r\leq s <t$ and $x_0$ in $\R^n$, there is a  $(Q^{\sigma, \lambda}_{r,x_0})$-null set $N$ such that for all $\omega$ in $N^c$, $E_{Q^{\sigma, \lambda}_{r,x_0}}(\phi(X_t)|{\cal B}^r_{s})(\omega)= U^{\sigma, \lambda}_t(\phi)(s,X_{s}(\omega))$. For given $a$ and $\lambda$, the modulus of continuity of 
$ U^{\sigma, \lambda}_t(\phi)$  at a point $(s,x)$ with $s<t$ depends only on the uniform bound of $\phi$.\\ 
\label{lemmacont}
\end{proposition}
{\bf Proof}
From \cite{SV2} Theorem 7.1, given $s<t$,   for every $\phi$ borelian bounded on $\R^n$, there is a function $U^{\sigma, \lambda}_{t}(\phi)$ such that for all $x$ in $R^n$, 
$E_{Q^{\sigma, \lambda}_{s,x}}(\phi(X_t))=U^{\sigma, \lambda}_{t}(\phi)(s,x)$. The map $(s,x)\in [0,t[\times \R^n \rightarrow U^{\sigma, \lambda}_{t}(\phi)(x)$ is continuous bounded and for given $a$ and $\lambda$, the modulus of continuity of 
$ U^{\sigma, \lambda}_{t}(\phi)$  at a point $(s,x)$ with $s<t$ depends only on the uniform bound on $\phi$.\\ 
It follows then from Lemma \ref{lemmapc1} that for every $r \leq s <t$, 
forall $\omega$ in $N^c$, $E_{Q^{\sigma, \lambda}_{r,x_0}}(\phi(X_t)|{\cal B}^r_{s})(\omega)=  U^{\sigma, \lambda}_t(\phi)(s,X_{s}(\omega))$ 
$\square $

\begin{proposition}Let $\sigma$ be continuous bounded and $a=\sigma \sigma^*$ strictly elliptic. Let  $\lambda$ be borelian bounded. 
\begin{enumerate}
\item For all sequence $s_n<s$ with limit $s$, the sequence  $Q^{\sigma, \lambda}_{s_n,x_n} (\pi^{s_n,s})^{-1}$ admits the limit  $Q^{\sigma, \lambda}_{s,x}$ for the weak topology. 
\item For $f$ continuous bounded on $\R^n$,  the map $U^{\sigma, \lambda}_t(f)(s,x)$, already defined in Proposition \ref{lemmacont} on $[0,t[ \times \R^n$, has a continuous extension to $[0,t] \times \R^n$ such that  for all $x$, $U^{\sigma, \lambda}_t(f)(t,x)=f(x)$.
\end{enumerate}
\label{propcontfell0}
\end{proposition}
{\bf Proof}
\begin{enumerate}
\item The proof is similar to that of Proposition \ref{lemmawrc}
\item Let $f$ be continuous bounded on $\R^n$. Let $t>0$. The map  $\omega \rightarrow f(X_t(\omega))$ is continuous bounded. Thus  from the weak convergence proved in $1.$,  it follows that for all  sequence $(s_n,x_n) \in [0,t[ \times \R^n$, with limit $(s,x)$,  the sequence $E_{Q^{\sigma \lambda}_{s_n,x_n}}(f(X_t)$ has the limit  $f(x)$. The result follows then from  Proposition \ref{lemmacont}. $\square$ \end{enumerate}


\subsection{Stable set of probability measures in Markovian setting}
\label{Stable}
From now on  $\sigma(t,x)$ is a given continuous bounded function, $a=\sigma \sigma^*$  is strictly elliptic. Uniform ellipticity is not assumed.  For given $r \geq 0$ and $y$ in $\R^n$,  we will consider probability measures all equivalent  to the probability measure ${Q}^{a}_{r,y}$. 
Recall that from \cite{BN03} the main ingredient in order to   construct a sublinear dynamic procedure time consistent  for deterministic times in $L^{\infty}(\Omega,{\cal B}^r,{Q}^{a}_{r,y})$ is  a set of equivalent probability measures stable by composition and stable by bifurcation   (cf \cite{BN03} definition 4.1) The definition is recalled in the appendix  ( definition  \ref{defstable}). Notice that a  more constraining condition is the m-stability introduced by Delbaen \cite{D}. It corresponds to time consistency for stopping times.\\
For the definition of multivalued mapping and of predictable multivalued mapping, we refer to \cite{D}. We introduce  now a definition for a  multivalued Borel mapping. 
\begin{definition}
A   multivalued Borel  mapping $\Lambda$ from $\R_+ \times \R^n$ into  a topological  vector space $E$ is a map $\Lambda$ defined on $\R_+ \times \R^n$ with values in the set of  subsets of  $E$ such that the graph of $\Lambda$ i.e. $$\{(t,x, \mu)|\mu \in \Lambda (t,x)\}$$
belongs to the Borel  $\sigma$-algebra ${\cal B}(\R_+ \times \R^n \times E)$, and such that $0 \in \Lambda(t,x)$, for all $(t,x)$.\\
It can have additional properties:
\begin{enumerate}
\item $\Lambda$ is convex if 
$\forall (t,x) \in  \R_+ \times \R^n$, $\Lambda(t,x)$ is a convex subset of $E$.
\item $\Lambda$ is closed if for all $(t,x)$, $\Lambda(t,x)$ is closed.
\end{enumerate}
\label{defBor}
\end{definition}
 Given $\sigma$, given $r \geq 0$ and $y$ in $\R^n$,   we want  to associate to a multivalued Borel mapping $\Lambda$  a set ${\cal Q}_{r,y}(\Lambda)$ of probability measures all equivalent with the probability measure $Q^a_{r,y}$ on ${\cal B}^r_t$ such that the process $(X_s)$ is a  strong Feller process with respect to every probability measure in ${\cal Q}_{r,y}(\Lambda)$. Therefore we introduce the following set:
\begin{definition} 
Let $\sigma$ be continuous bounded and $a=\sigma \sigma^*$ strictly elliptic. Let  $\Lambda$  be a multivalued Borel mapping  from $\R_+ \times \R^n$ into $\R^n$.
\begin{itemize}
\item
Let  $L(\Lambda)$ be the set of  bounded Borelian maps $\lambda: \R_+ \times \R^n \rightarrow \R^n$ such that for every $(t,x)$, $\lambda(t,x) \in \Lambda(t,x)$. 
 \item  For given $y$ in $\R^n$, and $r \geq 0$,  define the set of probability measures  ${\cal Q}_{r,y}(\Lambda)$
\begin{equation}
{\cal Q}_{r,y}(\Lambda)  =  \{Q^{\sigma, \lambda}_{r,y}, \lambda \in L(\Lambda)\}
\label{eqprob0}
\end{equation}
Recall that for $T>r$,
\begin{equation}
(\frac{dQ^{\sigma, \lambda}_{r,y}}{dQ^{a}_{r,y}})_{{\cal B}^r_T}  = exp(\int_r^T \lambda^*(t,X_t)\sigma^{-1}(t,X_t)dX_t-\frac{1}{2} \int_r^T
||\lambda(t,X_t)||^2dt)
\label{eqprob}
\end{equation}
\item For given $y$ in $\R^n$, denote also 
\begin{equation}
\tilde {\cal Q}_{r,y}(\Lambda)  =  \{Q^{\sigma, \lambda}_{r,y},\; \lambda \; \text{continuous},\;\lambda \in L(\Lambda)\}
\label{eqprob01}
\end{equation}
\end{itemize}
\label{def3}
\end{definition}
In order to construct a time consistent dynamic process for deterministic times, we need a  set of probability measures stable by composition and stable by bifurcation (cf \cite{BN03} or the appendix, definition  \ref{defstable}). 
Notice that the above set ${\cal Q}_{r,y}(\Lambda)$ is not stable by bifurcation. Indeed let  $\lambda$ and $\beta$ belonging to $L(\Lambda)$. Let $r \leq s<t$. Let  $B$ be  a Borelian set in $\R^n$, thus $A =\{\omega \in \Omega \;| \; X_s(\omega) \in B\}$ is  ${\cal B}^r_s$ measurable. $1_A(\omega)$ is not a function of $X_t(\omega)$ for $t>s$. Thus in general there is no  Borelian map $\alpha$ such that the process  defined as $\lambda(t,X_t(\omega))1_A+
\beta(t,X_t(\omega))1_{A^c}$ for $t > s$ and $\lambda(t,X_t(\omega)))$ for $t \leq  s$ can be written as $\alpha(t,X_t(\omega))$.  We need  to construct a stable set containing ${\cal Q}_{r,y}(\Lambda)$. We will define it as a subset of a certain stable set. For every $(t, \omega) \in \R_+ \times \Omega$, denote $\tilde \Lambda (t,\omega)=\Lambda(t,X_t(\omega))$.  
 
\begin{definition} 
Let $\sigma$ be continuous bounded and $a=\sigma \sigma^*$ strictly elliptic. For every bounded  process $\mu$ ${\cal B}^r_T$ predictable,  denote $P^{\sigma,  \mu}_{r,y}$ the probability measure equivalent with $Q^{a}_{r,y}$, on ${\cal B}^r_T$ with  Radon Nikodym derivative  
\begin{equation}
(\frac{dP^{\sigma, \mu}_{r,y}}{dQ^{a}_{r,y}})_{{\cal B}^r_T}=exp(\int_r^T \mu^*(u,\omega)\sigma^{-1}(u,X_u)(\omega)dX_u(\omega)-\frac{1}{2} \int_r^T
||\mu(u,\omega)||^2du)
\label{eqPmu}
\end{equation}
\label{defPmu}
\end{definition}

\begin{definition} .
Let $0\leq r \leq T$. Let  $\mu$ be a  bounded process ${\cal B}^r_T$-measurable. One says that   $\mu$ takes values in $\tilde \Lambda$ if  for all $r \leq s \leq T$, for all $\omega$, $\mu(s,\omega) \in \tilde \Lambda(s, \omega)$. Denote ${\cal M}_{r,y}(\tilde \Lambda)$  the set of all probability measures $P^{\sigma, \mu}_{r,y}$, on ${\cal B}^r_T$  where $\mu$ is ${\cal B}^r_T$ predictable and  takes  values in $\tilde \Lambda$. 
\label{def3bis}
\end{definition}
\begin{lemma}
 The set ${\cal M}_{r,y}(\tilde \Lambda)$ is a stable  set of probability measures on $(\Omega^r, {\cal B}^r_T)$. For all $k>0$, 
\begin{equation}
{\cal M}_{r,y}(\tilde \Lambda)^k= \{P^{\sigma, \mu}_{r,y} \in {\cal M}_{r,y}(\tilde \Lambda)\;|\;\;||\mu|| \leq k\}
\label{eqborne}
\end{equation}
 is also stable.\\
Furthermore if the mulivalued Borel mapping $\Lambda$ is convex,  ${\cal M}_{r,y}(\tilde \Lambda)$ and  ${\cal M}_{r,y}(\tilde \Lambda)^k$ are convex.
\label{lemmastable}
\end{lemma}
{\bf Proof}
 As already noticed, for given $ T \geq r$, $P^{\sigma, \mu}_{r,y}$ is a probability measure on  ${\cal B}^r_T$  equivalent with  $Q^a_{r,y}$.  The stability by composition and bifurcation  is easily verified. The convexity (in case $\Lambda$ is convex) follows from   the proof of Theorem 3  of \cite{D}.  The properties of  ${\cal M}_{r,y}(\tilde \Lambda)^k$ result then from Lemma 4 and Lemma 5 of \cite{D}.
\hfill $\square $\\
\begin{lemma} Let $0 \leq r \leq T$. The restriction to $(\Omega^r,  {\cal B}^r_T)$  of   ${\cal Q}_{r,y}(\Lambda)$ (resp. $\tilde {\cal Q}_{r,y}(\Lambda)$), given by (\ref{eqprob0}) (resp (\ref{eqprob01})),  is a set of equivalent probability measures on $(\Omega^r, {\cal B}^r_T)$.   There is a minimal stable set of probability measures containing ${\cal Q}_{r,y}(\Lambda)$ (resp. $\tilde {\cal Q}_{r,y}(\Lambda)$). We denote it $({\cal Q}_{r,y})_S(\Lambda)$ (resp. $(\tilde {\cal Q}_{r,y})_S(\Lambda)$).$(\tilde {\cal Q}_{r,y})_S(\Lambda)$) is a subset of $({\cal Q}_{r,y})_S(\Lambda)$. Every probability measure in $({\cal Q}_{r,y})_S(\Lambda)$  is equal to $P^{\sigma, \mu}_{r,y}$ for a certain  bounded  predictable process $\mu$ $\tilde \Lambda$ valued.
\label{lemmastab}
\end{lemma}
{\bf Proof}  Every probability measure in ${\cal Q}_{r,y}(\Lambda)$ belongs to ${\cal M}^{a}_{r,y}(\tilde \Lambda)$. Furthermore from 
Lemma \ref{lemmastable},  ${\cal M}^{a}_{r,y}(\tilde \Lambda)$ is stable.  The intersection of all stable sets of probability measures containing ${\cal Q}_{r,y}(\Lambda)$ is stable.  It is the minimal stable set of probability measures containing ${\cal Q}_{r,y}(\Lambda)$. It is a subset of ${\cal M}^{a}_{r,y}(\tilde \Lambda)$.
\hfill $\square $\\
 Notice that the set $({\cal Q}_{r,y})_S(\Lambda)$ is not closed for the weak topology in general.\\
We describe now the elements of  $({\cal Q}_{r,y})_S(\Lambda)$.
\begin{definition}
Let $0 \leq r \leq T$. ${\cal S}^r_T(\Lambda)$ (resp ${\tilde{\cal S}}^r_T(\Lambda)$) denotes the set of ${\cal B}^r_T$ measurable processes $\mu$ such that:\\ 
There is a finite subdivision $r=s_0<s_1<...<s_n=T$.\\  For all $i \in \{0,1,...n-1\}$ there is a finite set $I_i$, a finite partition $(A_{i,j})_{j \in I_i}$ of $\Omega$ into  ${\cal B}^r_{s_i}$-measurable sets, and Borelian (resp. continuous) bounded maps $\lambda_{i,j}$ in $L(\Lambda)$ such that 
\begin{equation}
\forall s_i  < u  \leq  s_{i+1}, \;\forall \omega \in \Omega,\;\mu(u,\omega)=\sum_{j \in I_i} \lambda_{i,j}(u,X_u(\omega))1_{A_{i,j}}(\omega)
\label{eqstable}
\end{equation}
\label{Slambda}
\end{definition}
\begin{proposition}
Given $0 \leq r \leq T$, the set $({\cal Q}_{r,y})_S(\Lambda)$ (resp. $(\tilde {\cal Q}_{r,y})_S(\Lambda)$),  is the set of all probability measures  $P^{\sigma, \mu}_{r,y}$  for some process $\mu$ belonging to ${\cal S}^r_T(\Lambda)$ (resp ${\tilde{\cal S}}^r_T(\Lambda)$).
\label{propdes}
\end{proposition}
{\bf Proof}
It is enough to do the proof for $({\cal Q}_{r,y})_S(\Lambda)$.\\Let $\mu$ satisfying (\ref{eqstable}).  Notice that for all $X$ ${\cal B}^r_{s_{i+1}}$ measurable, 
\begin{equation}
E_{P^{\sigma, \mu}_{r,y}}(X|{\cal B}^r_{s_i})=\sum_{j \in I_i} 1_{A_{i,j}}E_{Q^{\sigma, \lambda_{i,j}}_{r,y}}(X|{\cal B}^r_{s_i})
\label{eqdec}
\end{equation}
Thus it follows by induction using  the stability property (cf  appendix, definition  \ref{defstable})  that every  $P^{\sigma, \mu}_{r,y}$ where $\mu$ satisfies (\ref{eqstable})  belongs to the stable set $({\cal Q}_{r,y})_S(\Lambda)$. On the other hand it is easy to verify that the set $\{ P^{\sigma, \mu}_{r,y}:\;\mu \in {\cal S}^r_T(\Lambda)\}$ is stable.
\hfill $\square $\\

\begin{lemma}
Every set ${\cal Q}$ of equivalent probability measures on ${\cal B}^r_T$ stable by composition and bifurcation for deterministic times is stable by composition for stopping times taking a finite number of real  values.
\label{lemmastopfin}
\end{lemma}
{\bf Proof}
 Let $\sigma$ be a stopping time taking a finite number of values. Thus $\sigma$ can be written  $\sigma= \sum_{i=1}^n s_i 1_{A_i}$, where $A_i$ is a partition of $\Omega$,  $r \leq s_1 < s_2 <... s_n \leq T$, $A_i$ is ${\cal B}^r_{s_i}$ measurable.. Let $R$ and $Q$ in ${\cal Q}$. Denote $S_{\sigma}$ the probability measure in ${\cal Q}$ with Radon Nykodym derivative $\frac{dS_{\sigma}}{dP}=\frac{\frac{dQ}{dP}}{(\frac{dQ}{dP})_{\sigma}}(\frac{dR}{dP})_{\sigma}$. From the stability by composition and the stability by bifurcation for deterministic times, it follows easily by iteration  that $S_{\sigma}$ belongs to ${\cal Q}$.  \hfill $\square $\\

\subsection{Time-consistent sublinear  procedures}
\begin{remark}
$X_t$ being a strong Feller process for the probability measure $Q^{a}_{r,y}$, it follows from \cite{Pro} that the $Q^{a}_{r,y}$-completed filtration is right-continuous. Every  set of probability measures (equivalent to $Q^{a}_{r,y}$ on ${\cal B}^r_s$) stable for the filtration ${\cal B}^r_s$ is also stable for the completed filtration. 
In the following $({\cal B}^r_s)_{r \leq s}$  denotes the $Q^{a}_{r,y}$ completed filtration.
\end{remark}

\begin{proposition}
Let $0 \leq r \leq T$. Let  $({\cal Q}_{r,y})_S(\Lambda)$   be the stable set of probability measures generated by   $({\cal Q}_{r,y})(\Lambda)$ as in Lemma \ref{lemmastab}. Let $r \leq s \leq t \leq T$, $s$ and $t$  being  $r$-stopping times taking a finite number of  values. For all   $Y$ 
 in $L^{\infty}(\Omega^r,{\cal B}^r_{t},Q^{a}_{r,y})$, the formula 
\begin{equation}
\Pi^{r,y}_{s,t}(Y)=\esssup_{P^{\sigma, \mu}_{r,y} \in ({\cal Q}_{r,y})_S(\Lambda) } E_{P^{\sigma, \mu}_{r,y}}(Y|{\cal B}^r_{s})
\label{eqtcd0}
\end{equation}
defines  a dynamic process  time consistent for stopping times taking a finite number of real values.
Furthermore for given   $s$ and $t$, the map $\Pi^{r,y}_{s,t}$ defined on $L^{\infty}(\Omega^r,{\cal B}^r_{t},Q^{a}_{r,y})$ is sublinear  monotone normalized continuous from below.\\
For given  $0 \leq r \leq t$ and $Y$  in  $L^{\infty}(\Omega^r,{\cal B}^r_{t},Q^{a}_{r,y})$, the process $(\Pi^{r,y}_{s,t}(Y))_{r \leq s \leq t}$ is a $Q^{a}_{r,y}$ supermartingale and admits a c\`adl\`ag version. 
\label{prop1-0}
\end{proposition}

{\bf Proof} The set $({\cal Q}_{r,y})_S(\Lambda)$ being stable, the first part of the statement follows from Theorem 4.4   of \cite{BN03}.\\
The  proof of the regularity of paths which was given  in \cite{BN04} Theorem 3 for normalized dynamic processes time consistent for stopping times can be  extended to normalized convex (and thus to sublinear)  processes which are time consistent  for stopping times taking a finite number of real  values.   
\begin{remark} In the preceding proposition, the set $({\cal Q}_{r,y})_S(\Lambda)$ can be replaced by $(\tilde{\cal Q}_{r,y})_S(\Lambda)$. The result remains true.
\label{rem2}
\end{remark}

\begin{corollary}
The definition of $\Pi^{r,y}_{s,t}(Y)$ can be extended to random variables $Y$ ${\cal B}^r_t$-measurable which are only $Q^a_{r,y}$-essentially bounded from below:\\ $\Pi^{r,y}_{s,t}(Y)=\lim_{n \rightarrow \infty} \Pi^{r,y}_{s,t}(Y\wedge n)$.
For every such $Y$, $\Pi^{r,y}_{s,t}(Y)$ satisfies (\ref{eqtcd0}), and the process $(\Pi^{r,y}_{s,t}(Y)_s)$ is optional.
\label{corop}
\end{corollary}
{\bf Proof}
Let  $Y$ be ${\cal B}^r_t$-measurable and $Q^a_{r,y}$-essentially bounded from below, $Y$ is the increasing limit of $Y_n=Y \wedge n$ as $n$ tends to $\infty$.  $\Pi^{r,y}_{s,t}(Y)$ is then defined as the increasing limit of $\Pi^{r,y}_{s,t}(Y_n)$.  As we already know that for given $s$ and $t$, $\Pi^{r,y}_{s,t}$ defined on bounded random variables by formula (\ref{eqtcd0}) is continuous from below, the extended definition coincides with the previous one on essentially bounded random variables.\\
From the stability property of  $({\cal Q}_{r,y})_S(\Lambda)$, the set $\{E_Q(Y_n|{\cal B}^r_s)\}$ is a lattice upward directed. Thus   $\Pi^{r,y}_{s,t}(Y_n)$ is the increasing limit of a sequence $E_{Q_{n,k}}(Y_n)$ as $k \rightarrow \infty$. It is then easy to see  that $\Pi^{r,y}_{s,t}(Y)=\sup_{n,k} E_{Q_{n,k}}(Y)$
and to deduce that $\Pi^{r,y}_{s,t}(Y)$ satisfies (\ref{eqtcd0}). From  Proposition \ref{prop1-0} for every $n$ one can choose a c\`adl\`ag version of the process $\Pi^{r,y}_{s,t}(Y_n)_s$. Thus the map  $(s,\omega) \rightarrow \Pi^{r,y}_{s,t}(Y_n)(\omega)$ is measurable for the optional $\sigma$-algebra. It follows that  $s \rightarrow \Pi^{r,y}_{s,t}(Y)$ is an optional  process.
 \hfill $\square $ 

\section{Time-consistent convex   procedures  in  Markovian framework} 
\label{sectionTCC}
In Section \ref{secFellpro} we have constructed time consistent sublinear procedures  associated to a stable set of probability measures  in Markovian setting.  In this Section,   we want to construct time consistent convex dynamic procedures in Markovian settting.
The map $\sigma$  is given continuous bounded  and $a=\sigma \sigma^*$ is strictly elliptic. Uniform ellipticity is not assumed.
We assume also that a closed convex multivalued Borel mapping $\Lambda$ is given. 
 We have  constructed   in Section {\ref{Stable}  a stable set of probability measures which are all of the form $P^{\sigma, \mu}_{r,y}$ where $\mu$ is a bounded  process ${\cal B}^r_T$ measurable taking values in $\tilde \Lambda$ (cf Definition(\ref{defPmu}) for the definition of $P^{\sigma, \mu}_{r,y}$). 
We want now to construct families of   penalties  $\alpha_{st}(P^{\sigma, \mu}_{r,y})$ for $r \leq s \leq t \leq T$. The penalties will depend not only on the Borel mapping $\Lambda$ but also on a Borel measurable function $g$ with domain $\Lambda$.\\
Let  $g: \R^+ \times \R^n \times \R^n  \rightarrow \R \cup \{+\infty\}$ be  a Borel-measurable function such that   for all $(t,x) \in \R^+ \times \R^n$, 
$\Lambda(t,x) \subset\{y\in \R^n |\; g(t,x,y)<\infty\}$.  

 Define $f$ as follows: 
\begin{equation}
\forall z \in \R^d \;\;f(t,x,z)=\sup_{\lambda \in \Lambda(t,x) } (-z.\lambda-g(t,x,\lambda))
\label{eqdual}
\end{equation} 
The following lemma is straightforward:
\begin{lemma}
For all $(t,x)$, 
$f(t,x,.)$ is a  closed convex function  which is the dual transform of the function $\overline g (t,x,.)$ where 
\begin{eqnarray}
\overline g(t,x,\lambda)& = & g(t,x,\lambda)\;\; if \; \lambda \in \Lambda(t,x) \nonumber \\
                & = & +\infty \; else
\end{eqnarray}
For every $(t,x)$ $dom(\overline g(t,x,.))=\Lambda(t,x)$\\
If $g(t,x,0)=0 \;\forall (t,x)$, $f$ takes values  in $[0,\infty]$.\\ 
If $\overline g$ takes values in $[0,\infty]$ and satisfies $\forall (t,x),\; inf_{\lambda \in \Lambda(t,x)} g(t,x,\lambda)=0$ then for all $(t,x)$, $f(t,x,0)=0$.
\label{lemmadual}
\end{lemma}

Notice that, since $\Lambda$ is a closed convex multivalued Borel mapping, replacing  $g$ by $\overline g$, one can always assume that  for all  $(t,x)$, $dom(g(t,x,.)=\{\lambda \in \R^d\;| g(t,x,\lambda)<\infty \}$ is closed, convex and equal to $\Lambda(t,x)$.\\
In the following section we construct penalties from the above function $g$.


\subsection{Penalties }
\label{secpen}
In all the following, $\Lambda$ is a closed convex multivalued Borel mapping, and  $g: \R^+ \times \R^n \times \R^n \rightarrow  \R \cup \{\infty\}$ is a Borelian  map  such that  for all $(t,x), \;dom g(t,x,.)=\{y,\;g(t,x,y)<\infty\}=\Lambda(t,x)$.

\begin{definition} 
$g$ has polynomial growth on $\Lambda$ if there is $C>0$ and $m \in \N$ such that 
\begin{equation}
\sup_{y \in \Lambda(u,x)} |g(u,x,y)| \leq C(1+||x||^m)
\label{eqpg}
\end{equation}
\label{defpg}
\end{definition}

\begin{definition}
Assume that $g$ is non negative or has polynomial growth on $\Lambda$. Let $0 \leq r \leq T $. For all bounded  process $\mu$  $\tilde \Lambda$-valued with left continuous paths admitting right limits (c\`agl\`ad),  for all  
$r$-stopping times $r \leq s  \leq t \leq T$, taking a finite number of values  define the penalty  $\alpha_{s,t}(P^{\sigma, \mu}_{r,y})$ as follows
\begin{equation}
\alpha_{s,t}(P^{\sigma, \mu}_{r,y})=E_{P^{\sigma, \mu}_{r,y}}( \int\limits_{s}^{t} g(u,X_u(\omega),\mu(u,\omega)) d u|{\cal B}^r_{s})
\label{eqP1_}
\end{equation}
\label{defpen0}
\end{definition}
We have introduced in  \cite{BN03}, definition 4.3  the definition of local property and of cocycle condition for the penalty.  These definitions are recalled in the appendix (Definition \ref{defpen}).  
\begin{proposition}
\begin{itemize}
\item
\begin{enumerate}
\item In case  $g$ has polynomial growth on $\Lambda$, equation (\ref{eqP1_}) defines an element of $L^1(P^{\sigma, \mu}_{r,y})$ and of $L^1(Q^{a}_{r,y})$ for all bounded $\mu$.
\item In case  $g$ is non negative,  equation (\ref{eqP1_}) defines a non negative ${\cal B}^r_s$ random variable.
\end{enumerate}
\item 
The  penalty defined in (\ref{eqP1_})  satisfies the cocycle condition for every $P^{\sigma, \mu}_{r,y}$: Let $s_0$, $s$ and $t$ be $r$-stopping times taking a finite number of real values, $r \leq s_0 \leq s \leq t$
\begin{equation}
 \;\;\alpha_{s_0,t}(P^{\sigma, \mu}_{r,y})=\alpha_{s_0,s}(P^{\sigma, \mu}_{r,y})+E_{P^{\sigma, \mu}_{r,y}}(\alpha_{s,t}(P^{\sigma, \mu}_{r,y})|{\cal B}^{r}_{s_0})
\label{eqcoc}
\end{equation}
\item The penalty defined in (\ref{eqP1_}) is local. 
\item If $g(t,x,0)=0\; \forall (t,x) \in \R^+ \times \R^n$, 
The probability measure $Q^{a}_{r,y}=P^{\sigma,0}_{r,y}$ has zero penalty. 
\end{itemize}
\label{lemmapen}
\end{proposition}
{\bf Proof}
\begin{itemize}
\item 
\begin{enumerate}
\item
Assume that the function $g$ has polynomial growth on $\Lambda$ with growth exponent $m$. Choose $q \geq 1$ such that $mq \geq 2$. Let $p$ be the conjugate exponent of $q$.
The process $\mu$ being bounded,  
$(\frac{dP^{\sigma, \mu}_{r,y}}{dQ^a_{r,y}})_{{\cal B}^r_t}$ belongs to $L^p(Q^a_{r,y})$. Thus from H\"older inequality, 
\begin{equation}
E_{P^{\sigma, \mu}_{r,y}}(\sup_{s \leq u \leq t}||X_t||^m) \leq K'E_{Q^a_{r,y}}(\sup_{s \leq u \leq t}||X_t||^{mq})^{\frac{1}{q}}
\label{eqHold}
\end{equation}
 It follows from Proposition \ref{propKryl}  that   the penalty $\alpha_{s,t}(P^{\sigma, \mu}_{r,x})$ is always well defined and belongs to   $L^1(Q^{a}_{r,y})$ and  $L^1(P^{\sigma, \mu}_{r,y})$.
\item The case $g$ non negative is trivial.
\end{enumerate}
\item
 The cocycle condition (\ref{eqcoc}) follows easily from the definition (\ref{eqP1_}).
\item We prove now  that the penalty $\alpha$ is local. The probability measures  $P^{\sigma, \mu}_{r,y}$ and $P^{\sigma, \nu}_{r,y}$ are equivalent  to   $Q^{a}_{r,y}$ on ${\cal B}^r_t$. 
Let $A$ be  ${\cal B}^r_{s}$-measurable. Assume that  $\forall X \in L^{\infty}({\cal B}^r_{t}),\; E_{P^{\sigma,\nu}_{r,y}}(X|{\cal B}^r_{s})1_A= E_{P^{\sigma, \mu}_{r,y}}(X|{\cal B}^r_{s})1_A$. The Radon Nikodym derivatives of $P^{\sigma, \mu}_{r,y}$  and $P^{\sigma,\nu}_{r,y}$ on ${\cal B}^r_t$ are given by equation (\ref{eqPmu}). It follows that $Q^a_{r,y}$ a.s. 
$$1_A(\int_s^t \mu^*(u,\omega)\sigma^{-1}(u,X_u)(\omega)dX_u(\omega)-\frac{1}{2} \int_s^t
||\mu(u,\omega)||^2du)=$$
$$1_A(\int_s^t \nu^*(u,\omega)\sigma^{-1}(u,X_u)(\omega)dX_u(\omega)-\frac{1}{2} \int_s^t
||\nu(u,\omega)||^2du)$$
$\mu$ and $\nu$ being c\`agl\`ad processes, it follows that for almost all $\omega$ in $A$,  $ \mu(u, \omega)=\nu(u,\omega)$ for all $s \leq u <t$.
 From (\ref{eqP1_}) it follows that $\alpha_{s,t}(P^{\sigma, \mu}_{r,y})1_A=\alpha_{s,t}(P^{\sigma, \nu}_{r,y})1_A$. Thus the penalty $\alpha$ is local.
\item The last point follows easily from the definition of the penalty.  \hfill $\square $ 
\end{itemize}


\subsection{Normalized time consistent convex  procedure associated to a multi-valued Borel mapping and a  non negative Borel map}

\begin{proposition}
Let  $(\tilde {\cal Q}_{r,y})_S({\Lambda})$  be  as in Lemma \ref{lemmastab}. Assume that $g$ is non negative, and that for all $(u,x)$, $g(u,x,0)=0$. Let $s$ and  $t$ be $r$-stopping times taking a finite number of  values $r \leq s \leq t$. The formula 
\begin{equation}
\Pi^{r,y}_{s,t}(Y)=\esssup_{P^{\sigma, \mu}_{r,y} \in (\tilde {\cal Q}_{r,y})_S({\Lambda})} (E_{P^{\sigma, \mu}_{r,y}}(Y|{\cal B}^r_{s})- \alpha_{s,t}(P^{\sigma, \mu}_{r,y}))
\label{eqtcd000}
\end{equation}
where $\alpha_{s,t}(P^{\sigma, \mu}_{r,y})$ is given by equation (\ref{eqP1_})
defines  a normalized convex dynamic process on  $ L^{\infty}(\Omega,({\cal B}^r_{t})^+,Q^{a}_{r,y})$ time consistent for $r$-stopping times taking a finite number of values.
Furthermore for fixed $s$ and $t$, the map $\Pi^{r,y}_{s,t}$ defined on  $ L^{\infty}(\Omega,({\cal B}^r_{t})^+,Q^{a}_{r,y})$ is convex  monotone  continuous from below.\\
\label{prop1}
\end{proposition}
{\bf Proof} Notice that for all bounded $Y$, 
$$-||Y||_{\infty} \leq E_{Q^{a}_{r,y}}(Y|{\cal B}^r_{s}) \leq \Pi^{r,y}_{s,t}(Y)\leq \esssup_{P^{\sigma, \mu}_{r,y} \in (\tilde {\cal Q}_{r,y})_S({\Lambda})} E_{P^{\sigma, \mu}_{r,y}}(Y|{\cal B}^r_{s})$$
Thus for all $r  \leq s \leq t,\; ||\Pi^{r,y}_{s,t}(Y)||_{\infty} \leq ||Y||_{\infty}$.
The statement follows then from Theorem 4.4   of \cite{BN03}.  \hfill $\square $ \\
As in the sublinear case (cf Corollary(\ref{corop}) we have the following extension
\begin{corollary}
The definition of $\Pi^{r,y}_{s,t}(Y)$ can be extended to random variables $Y$ $({\cal B}^r_t)^+$-measurable which are only essentially bounded from below. \\ $\Pi_{s,t}(Y)=\lim_{n \rightarrow \infty} \Pi_{s,t}(Y\wedge n)$.
For every such $Y$, $\Pi_{s,t}(Y)$ satisfies (\ref{eqtcd000}).
\label{corop2}
\end{corollary}


\subsection{ General time consistent convex procedure associated to a multi-valued Borel mapping and a  Borel map}
In this section the function $g$ (and thus the penalty) is not assumed to be non negative. We assume now that the function $g(u,x,0)$ is bounded from above on $\R_+ \times \R^n$.
It follows that, given $T$, the penalty $\alpha_{s,t}(Q^{a}_{r,y})$ is bounded from above, uniformly  in $(s,t)$  such that   $0 \leq r \leq s \leq t \leq T$. 

\begin{proposition}
Assume that $\Lambda$ is a  closed convex multivalued Borel mapping. Assume that $g$ has polynomial growth and that $g(u,x,0)$ is bounded from above.  Let  $(\tilde {\cal Q}_{r,y})_S({\Lambda})$  be  as in Lemma \ref{lemmastab}. Given $r \leq s \leq t$, the formula 
\begin{equation}
\Pi^{r,y}_{s,t}(Y)=\esssup_{P^{\sigma, \mu}_{r,y} \in (\tilde {\cal Q}_{r,y})_S({\Lambda})} (E_{P^{\sigma, \mu}_{r,y}}(Y|{\cal B}^r_{s})- \alpha_{s,t}(P^{\sigma, \mu}_{r,y}))
\label{eqtcd}
\end{equation}
where $\alpha_{s,t}(P^{\sigma, \mu}_{r,y})$ is given by equation (\ref{eqP1_})
defines  a  convex   monotone map continuous from below on  the set of ${\cal B}^r_t$ measurable variables essentially bounded from below with values in  the set of ${\cal B}^r_s$ measurable variables essentially bounded from below. \\
Given $r$ and $y$, $\Pi^{r,y}_{s,t}$ is a convex  dynamic process time consistent for stopping times taking a finite number of values.
\label{prop1_1}
\end{proposition}
{\bf Proof}
For all $Y$ bounded from below, $(E_{Q^{a}_{r,y}}(Y|{\cal B}^r_{s})- \alpha_{s,t}(Q^{a}_{r,y}))$ is essentially bounded from below, thus it is the same for $\Pi^{r,y}_{s,t}(Y)$.

The result follows then from the proof of Theorem 4.4 of \cite{BN03} which can be adapted without  difficulty to the case of variable essentially bounded from below.
\hfill $\square $. 

\section{Strong Feller property of the time consistent convex dynamic procedure}
\label{secFeller}
The goal of this Section is to prove a Feller property for  the dynamic process $\Pi^{r,y}_{s,t}(h(X_t))$.   Recall that  from Section \ref{subsec1-3}, for all bounded $\lambda$, $r$ and $y$, the process $X_t$ is a Feller process with respect to the  probability measure $Q^{\sigma,\lambda}_{r,y}$.

\subsection{Feller property for the penalty associated to a Feller probability measure}
\begin{definition}
A real valued function $f$ defined on $\R^+ \times \R^n$ 
\begin{itemize}
\item
is a  Caratheodory function if it is Borelian and if for all $t$, $f(t,.)$ is continuous on $\R^n$.
\item  has polynomial growth on $I \times \R^n$ where $I$ is a subset of $\R^+$ if there is a constant $C>0$ and $m \in \N$ such that 
\begin{equation}
\forall (t,x) \in I \times \R^n,\;\;|f(t,x)| \leq C(1+||x||^m)
\label{eqpolg}
\end{equation}
\end{itemize}
\label{defCat}
\end{definition}
We prove a Feller property for Caratheodory functions $f$ with polynomial growth.
\begin{proposition} Let $\sigma$ be continuous bounded and $a=\sigma \sigma^*$ strictly elliptic. Let  $\lambda$ be borelian bounded.  Assume that $f$ is a   real valued   Caratheodory function with polynomial growth on  $[0,t] \times \R^n$. 
\begin{enumerate}
\item  There is a real valued continuous map $L(f)$   on $[0,t]\times \R^n$ such that 
\begin{equation}
E_{Q^{\sigma,\lambda}_{s,y}}(\int_s^t f(u, X_u) du )=L(f)(s,y) \;\;\;\forall s \in [0,t]\;\;and\; y \in \R^n
\label{eqFell47-3}
\end{equation}
\item For all $0 \leq r \leq s \leq  t$ and all $y \in \R^n$, there is a 
$Q^{\sigma, \lambda}_{r,y}$-null set N such that for all $\omega \in N^c$, 
\begin{equation}
E_{Q^{ \sigma, \lambda}_{r,y}}(\int_s^t f(u, X_u) du |{\cal B}^r_s)(\omega)=L(f)(s,X_s(\omega))
\label{eqFell47b-3}
\end{equation}
\end{enumerate}
Notice that for all $y \in \R^n$ $L(f)(t,y)=0$.
\label{propcont3}
\end{proposition}

{\bf Proof}
We prove statement {\it 1.} in three steps.
\begin{itemize}
\item
Let $0 \leq s \leq t$ and $y$ in $\R^n$. Define $L(f)(s,y)$ by equation (\ref{eqFell47-3}). The function $f$ being Borelian with  polynomial growth and $\lambda$ bounded it follows from Proposition \ref{propKryl}  that $L(f)(s,y)$ is a real number.
\item We prove the continuity of $L(f)$ at every point $(s,x)$ for $s<t$.\\
Let $(s,x) \in [0,t[ \times \R^n$.  Choose  $\eta>0$ such that $s+\eta<t$. By hypothesis $f$ has polynomial growth., thus $\forall (s',x') \in [0,s+\eta] \times \{y \in \R^n,\; ||y||\leq ||x||+1\}$, 
\begin{eqnarray}
E_{Q^{\sigma, \lambda}_{s',x'}}(\int_{s+\eta}^t|f(u,X_u)|1_{\{\sup_{s' \leq u \leq t}||X_u|| \geq k\}} du \leq \nonumber \\
E_{Q^{\sigma, \lambda}_{s',x'}}(\int_{s+\eta}^t C(1+||X_u||^m)(|1_{\{\sup_{s' \leq u \leq t}||X_u|| \geq k\}}du 
\label{eqcontf00}
\end{eqnarray}
Notice that $||X_u|| \geq k$ implies $||X_u-x'||\geq k-||x||-1$.  Let $\epsilon>0$. It follows then from Corollary \ref{corKry} that there is $k_0>0$ such that for $k \geq k_0$ and $(s',x') \in [0,s+\eta] \times \{y \in \R^n,\; ||y||\leq ||x||+1\}$,
\begin{equation}
E_{Q^{\sigma, \lambda}_{s',x'}}(\int_{s+\eta}^t|f(u,X_u)|1_{\{\sup_{s' \leq u \leq t}||X_u|| \geq k\}} du \leq \epsilon
\label{eqcontf0}
\end{equation}
Let $f_k=\sup(\inf(f,k),-k)$. From equation (\ref{eqpolg}), $f(u,X_u)(\omega) \neq f_k(u,X_u)(\omega)$ implies that $||X_u(\omega)||\geq (\frac{k}{C}-1)^{\frac{1}{m}}$. It follows from equation (\ref{eqcontf0}) that 
there is $k_1>0$ such that $\forall k \geq k_1$, $\forall (s',x') \in [0,s+\eta] \times \{y \in \R^n,\; ||y||\leq ||x||+1\}$, 
\begin{equation}
E_{Q^{\sigma, \lambda}_{s',x'}}(\int_{s+\eta}^t|f(u,X_u)-f_k(u,X_u)|)du \leq \epsilon.
\label{eqinte0}
\end{equation}
 $X_u$ being a continuous function of $\omega$ and $f(u,.)$ a continuous function on $\R^n$, it follows  that for all $k$, $\int_{s+\eta}^t f_k(u, X_u(\omega)) du $ is a continuous bounded function of $\omega$. \\
Let $(s_n, x_n)$ be a sequence with limit $(s,x)$. From Proposition (\ref{propcontfell0}), the sequence of probability measures $Q^{\sigma, \lambda}_{s_n,x_n}(\Pi^{s_n,s})^{-1}$ converges to  $Q^{\sigma, \lambda}_{s,x}$ for the weak topology. Thus there is $N>0$ such that for all $n \geq N$, 
\begin{equation}
|E_{Q^{\sigma, \lambda}_{s_n,x_n}}(\int_{s+\eta}^t f_{k_1}(u, X_u(\omega)) du - E_{Q^{\sigma, \lambda}_{s,x}}(\int_{s+\eta}^t f_{k_1}(u, X_u(\omega)) du |\leq \epsilon
\label{eqconv}
\end{equation}
It follows then easily from equations (\ref{eqinte0}) and (\ref{eqconv}) that 
\begin{equation}
\lim_{n \rightarrow \infty}E_{Q^{\sigma, \lambda}_{s_n,x_n}}(\int_{s+\eta}^t f(u, X_u(\omega)) du =E_{Q^{\sigma, \lambda}_{s,x}}(\int_{s+\eta}^t f(u, X_u(\omega)) du 
\label{eqconv2}
\end{equation}
From Proposition  \ref{propKryl}, there is a constant $K>0$ such that $\forall s' \leq t_1 \leq t$, and  $||x'|| \leq ||x||+1\}$, 
\begin{equation}
E_{Q^{\sigma, \lambda}_{s',x'}}(\int_{s'}^{t_1} |f(u, X_u(\omega))| du  \leq K(t_1-s')
\label{eqkcont}
\end{equation}
The continuity of $L(f)$ at $(s,x)$ follows from equations (\ref{eqconv2}) and (\ref{eqkcont}) applied with $t_1=s+\eta$ 
\item Continuity at $(t,x)$. Let  $(s_n,x_n)$ be a sequence with limit $(t,x)$.
 From equation (\ref{eqkcont})  (applied with $t_1=t$),  it follows that   $L(f)(s_n,x_n) \rightarrow 0=L(f)(t,x)$ as $n \rightarrow \infty$.
\end{itemize}
This proves {\it 1.}\\ Statement {\it 2.} is then a consequence of Lemma \ref{lemmapc1}.
\hfill $\square $

We introduce now the corresponding hypothesis on $g$.

\begin{definition}
Hypothesis $H_g$\\
\begin{enumerate}
\item 
$g:\R_+ \times \R^n \times \R^n \rightarrow \R$ is a ``Caratheodory function on $\Lambda$''\\
More precisely, $g$ is Borelian and for all $u$, the restrition of $g_u$ to $\{(x,y), y \in \Lambda(u,x)\}$ is continuous  ($g_u(x,y)=g(u,x,y)$).
\item $g$ has polynomial growth on $\Lambda$  (cf Definition \ref{defpg}).
\end{enumerate}
\label{defHg}
\end{definition}

\begin{corollary} Let $\sigma$ be continuous bounded and $a=\sigma \sigma^*$ strictly elliptic.   Let  $\lambda$ be a  bounded Caratheodory function $\Lambda$ valued.  Assume that $g$ satisfies Hypothesis $H_g$.
  For all $t>0$, there is a continuous  map $(s,x) \in [0,t]\times \R^n \rightarrow L^{\sigma,\lambda}_t(g)(s,x)$ such that  for all $0 \leq r \leq s \leq t$ and all $y \in \R^n$, there is a 
$Q^{\sigma, \lambda}_{r,y}$-null set N such that for all $\omega \in N^c$, 
 \begin{equation}
\alpha_{st}(Q^{\sigma, \lambda}_{r,y})=E_{Q^{\sigma, \lambda}_{r,y}}(\int_s^t g(u, X_u,\lambda(u,X_u)) du |{\cal B}^r_s)(\omega)=L^{\sigma,\lambda}_t(g)(s,X_s(\omega))
\label{eqFell47b}
\end{equation}
\label{colcont2b}
\end{corollary}
{\bf Proof}
The map $f(u,x)=g(u,x,\lambda(u,x))$ satisfies the hypothesis of Proposition \ref{propcont3}. This gives the existence of the function  $L^{\sigma,\lambda}_t(g)$ satisfying  the required conditions.

\subsection{Feller property for the dynamic convex procedure}

\begin{proposition} Let $\sigma$ be continuous bounded and $a=\sigma \sigma^*$ strictly elliptic.   Let  $\lambda$ be a  bounded Caratheodory function $\Lambda$ valued. Assume that the penalty is given by equation(\ref{eqP1_}) (Definition \ref{defpen0}) for some function $g$ satisfying hypothesis $H_g$ (Definition \ref{defHg}).
\begin{enumerate}
\item For all $\phi$ Borelian bounded on $\R^n$, and  $t>0$,  there is a   function $V^{\sigma,\lambda}_t(\phi)$ continuous  on $[0,t[\times \R^n$, bounded on $[0,t[\times \{x \in \R^n|\; ||x|| \leq K\}$ for all $K$, such that
\begin{equation}
E_{Q^{\sigma, \lambda}_{s,x}}(\phi(X_t))-\alpha_{st}(Q^{\sigma, \lambda}_{s,x})=V^{\sigma,\lambda}_t(\phi)(s,x)
\label{eqFell_0}
\end{equation}
For all $r \leq s \leq t$ and $y \in \R^n$,
\begin{equation}
E_{Q^{\sigma, \lambda}_{r,y}}(\phi(X_t)|{\cal B}^r_s)-\alpha_{st}(Q^{\sigma, \lambda}_{r,y})=V^{\sigma,\lambda}_t(\phi)(s,X_s)\;\; Q^{a}_{r,y} a.s.
\label{eqFell1}
\end{equation}
\item If $\phi$ is furthermore continuous on $\R^n$, $ V^{\sigma,\lambda}_t(\phi)$ is continuous  on $[0,t] \times \R^n$ with $ V^{\sigma,\lambda}_t(\phi)(t,x)=\phi(x)$
\end{enumerate}
\label{propcontfell}
\end{proposition}

{\bf Proof}
\begin{enumerate}
\item
The equality (\ref{eqFell_0}) follows from   Proposition \ref{lemmacont}  and Corollary \ref{colcont2b}, with 
$$V^{\sigma,\lambda}_t(\phi)=U^{\sigma, \lambda}_t(\phi)-L^{\sigma,\lambda}_t(g)$$ 
 $U^{\sigma, \lambda}_t(\phi)$  is continuous  on $[0,t[\times \R^n$ and bounded on $[0,t[\times \{x \in \R^n|\; ||x|| \leq K\}$ for all $K$. The function $L^{\sigma,\lambda}_t(g)$ is continuous  on $[0,t]\times \R^n$.\\
The equality (\ref{eqFell1}) is then a consequence of Lemma \ref{lemmapc1}.
\item  The result follows from the above proof and from Proposition \ref{propcontfell0} $\square$
\end{enumerate}
Following \cite{SV1} (cf Section 3) we say that a process $\mu$ is r-non-anticipating if for all $t \geq r$, $\mu(t,\omega)$ is ${\cal B}^r_t$-measurable. 
\begin{proposition}
Let  $P^{\sigma, \mu}_{r,y}$ be a probability measure on $(\Omega^r,{\cal B}^r_T)$ belonging to  $({\cal Q}_{r,y})_S(\Lambda)$. The process $\mu$ is bounded and r-non-anticipating.  The   probability measure $P^{\sigma, \mu}_{r,y}$ on $(\Omega^r,{\cal B}^r_T)$ is  the unique solution to the martingale problem
\begin{equation}
Y_t=exp \{(\theta)^*(X(t)-y)-\int_r^t \theta^*\sigma(u, X_u(\omega)) \mu(u,\omega)  du -\frac{1}{2} \int_r^t \theta^* a(u,X_u(\omega)) \theta du \} 
\label{eqmarmu}
\end{equation}
starting from $y$ at time $r$, i.e. such that $E_{P^{\sigma, \mu}_{r,y}}(\{X_r=y\})=1$. 
\label{uniqmar}
\end{proposition}
{\bf Proof}
We use the description of $({\cal Q}_{r,y})_S(\Lambda)$ given in Proposition \ref{propdes}. It follows that $\mu$ is bounded and r-non-anticipating. The probability measure $Q^{\sigma, \lambda_{i,j}}_{r,y}$ being solution to the martingale problem   
 (\ref{eqmar}) for $\lambda=\lambda_{ij}$ it follows by induction that $P^{\sigma \mu}_{r,y}$ is solution to the martingale problem (\ref{eqmarmu}).\\
The proof of the unicity of the solution to the martingale problem given in Section 6 (Theorem 6.2) of \cite{SV1} in the particular case where $\mu(u,\omega)=b(u, X_u(\omega))$ can be adapted without difficulties to the more general case  where $\mu$ is bounded and r-non-anticipating.  $\square$\\
From proposition \ref{propdes} and  the unicity result of Proposition \ref{uniqmar} we deduce as in Lemma \ref{lemmapc1} the following result
\begin{corollary}
Let $\nu_s$ in ${\cal S}^s_t(\Lambda)$ (cf Definition \ref{Slambda}).  Let $r<s$. Let $\tilde \nu_s$ in ${\cal S}^r_t(\Lambda)$   such that for all $u \geq s$ $\tilde \nu_s(u,\omega)=\nu_s(u,\omega)$. Let $(P^{\sigma,\tilde \nu_s}_{r,y})_{s,\omega}$ be  a regular conditional distribution of $P^{\sigma,\tilde \nu_s}_{r,y}$ given ${\cal B}^r_s$, then
\begin{equation}
(P^{\sigma,\tilde \nu_s}_{r,y})_{s,\omega}(\pi^{r,s})^{-1}=P^{\sigma, \nu_s}_{s,X_s(\omega)}\;\;\;Q^{a}_{r,y}\;a.s.
\label{eqrcd}
\end{equation}
\label{corrcd}
\end{corollary}

\begin{theorem}
 Let $r<t$. Let $\mu$ in $\tilde{\cal S}^r_t(\Lambda)$. 
For all $ h$  in  ${\cal B}_b(\R^n)$ (resp. ${\cal C}_b(\R^n))$, there is $\overline h$ continuous on  $([r,t[\times \R^n)$ bounded on $\{x, \; ||x|| \leq K\}$ (resp. $\overline h$ continuous on $[r,t]\times \R^n$  and $\overline h (t,x)=h(x)$) such that for all $r\leq s <t$ there is a process $\nu_s$ in $\tilde {\cal S}^s_t(\Lambda)$ satisfying  (\ref{eqin1}) \;and\; (\ref{eqin2})
\begin{equation}
\forall x, y  \in \R^n,\;\; 
E_{P^{\sigma, \nu_s}_{s,x}}(h(X_t))-\alpha_{s,t}(P^{\sigma, \nu_s}_{s,x})= \overline h(s,x)
\label{eqin1}
\end{equation}
\begin{equation}
 E_{P^{\sigma, \mu}_{r,y}}(h(t,X_t)|{\cal B}^r_s)-\alpha_{s,t}(P^{\sigma, \mu}_{r,y}) \leq  E_{P^{\sigma,  \tilde \nu_s}_{r,y}}(h(t,X_t)|{\cal B}^r_s)-\alpha_{s,t}(P^{\sigma,  \tilde \nu_s}_{r,y})=\overline h(s,X_s)
\label{eqin2}
\end{equation}
where $\tilde \nu_s(u,\omega)=0$ for $r \leq u \leq s$ and $\tilde \nu_s(u,\omega)=\nu_s(u,\omega)$ for $u>s$.
\label{thmMar}
\end{theorem}

{\bf Proof}
The proof is done in two steps. The first one is the construction of $\overline h$ given $\mu$.  The second one is the construction of $\nu_s$ given $\mu$ and $s$.
\begin{itemize}
\item First step: construction of $\overline h$.\\
 Let 
$r=s_0<s_1<...<s_n=t$ be the subdivision associated to $\mu$ as in Definition \ref{Slambda}. 
For  $u \in ]s_i,s_{i+1}]$,  $\mu(u,\omega)=\sum_{j \in I_i} 1_{A_{i,j}}(\omega)\lambda_{i,j}(u,X_u(\omega))$, where $\lambda_{ij}$ is continuous bounded. We define $\overline h$ recursively on $[s_i,s_{i+1}[$.
From Proposition \ref{propcontfell} for all $j \in I_{n-1}$ there is a  map $V^{\sigma,\lambda_{n-1,j}}_t(h)$ continuous on $[0,t[ \times \R^n$ such that equations (\ref{eqFell_0}) and (\ref{eqFell1}) are satisfied. let 
\begin{equation}
\overline h(s,x)=\sup_{j \in I_{n-1}}V^{\sigma,\lambda_{n-1,j}}_t(h)(s,x)\;\; \forall s \in [s_{n-1},t[
\label{hb1}
\end{equation}

$\overline h$ is continuous on $[s_{n-1},t[\times \R^n$.\\
Let $i+1<n$.
Assume  now that $\overline h$ has been defined as a continuous function on  $[s_{i+1},t[\times \R^n$. Let $\phi_i(x)=\overline h(s_{i+1},x)$.  From Proposition \ref{propcontfell}, for all $j \in I_{i}$ there is a  map $V^{\sigma,\lambda_{i,j}}_t(\phi_i)$ continuous on $[0,s_{i+1}] \times \R^n$ with $V^{\sigma,\lambda_{i,j}}_t(\phi_i)(s_{i+1},x)=\phi_i(x)=\overline h(s_{i+1},x)$. Let 
\begin{equation}
\overline h(s,x)=\sup_{j \in I_{i}}V^{\sigma,\lambda_{i,j}}_t(\phi_i)(s,x)\;\; \forall s \in [s_{i},s_{i+1}[
\label{hb2}
\end{equation}

$\overline h$ is continuous on $[s_{i},t[\times \R^n$. This ends the proof of the construction  of a continuous function $\overline h$
 associated to $\mu$ and $h$.
Notice that for all $s \in [s_{i},s_{i+1}[$, there is a partition of $\R^n$ in Borelian sets  $(C_{s,j})_{j \in I_i}$,  such that 
\begin{equation}
\overline h(s,x)=\sum_{j \in I_{i}}1_{C_{s,j}}(x)V^{\sigma,\lambda_{i,j}}_t(\phi_i)(s,x)\;\; 
\label{hb3}
\end{equation}

\item Second step: Given $s \in ]r,t]$, construction of the process $\nu_s$.\\ There is a unique $k$ such that $s \in ]s_k, s_{k+1}]$. For $s=r$, let $k=0$.
For $i>k$ for all $u \in ]s_i,s_{i+1}]$, define 
\begin{equation}
\nu_s(u,\omega)=\sum_{j \in I_{i}}1_{C_{s_i,j}}(X_{s_i}(\omega))\lambda_{i,j}(u,X_u(\omega))\;\
\label{eqnu1}
\end{equation}

And for  $u \in ]s,s_{k+1}]$, define 
\begin{equation}
\nu_s(u,\omega)=\sum_{j \in I_{k}}1_{C_{s,j}}(X_s(\omega))\lambda_{k,j}(u,X_u(\omega))\;\
\label{eqnu2}
\end{equation}

Let $\tilde \nu_s(u,\omega)=0$ for $u \leq s$.
For all $i>k$, the restriction of $\nu_s$ to $[s_i,t]$  belongs to $\tilde{\cal S}^{s_i}_t(\Lambda)$.   We still denote it $\nu_s$. From the construction of $\nu_s$, it follows recursively that for all $i>k$, and all $y \in \R^n$:
\begin{equation}
 E_{P^{\sigma, \nu_s}_{s_i,y}}(h(X_t))=\overline h(s_i,y)
\label{eqnu3b}
\end{equation}
The equality 
 \begin{equation}
E_{P^{\sigma, \tilde \nu_s}_{r,y}}(h(X_t)|{\cal B}^r_{s_{i}})-\alpha_{s_{i}t}(P^{\sigma, \tilde\nu_s}_{r,y})=\overline h(s_{i},X_{s_{i}})
\label{eqnu3}
\end{equation}
is deduced from equation (\ref{eqnu3b}) using Corollary \ref{corrcd}. \\
Endly
$E_{P^{\sigma, \nu_s}_{s,y}}(h(X_t)=\overline h(s,y)$ and  
\begin{equation}
E_{P^{\sigma, \tilde \nu_s}_{r,y}}(h(X_t)|{\cal B}^r_{s})-\alpha_{st}(P^{\sigma, \tilde \nu_s}_{r,y})=\overline h(s,X_{s})
\label{eqnu3c}
\end{equation}

Notice that from the expression of $\mu$ and from the definition of the penalty it follows that 
\begin{eqnarray}
E_{P^{\sigma,  \mu}_{r,y}}(h(X_t)|{\cal B}^r_{s_{n-1}})-\alpha_{s_{n-1}t}(P^{\sigma,  \mu}_{r,y})= \nonumber\\\sum_{j \in I_{n-1}}1_{A_{n-1,j}}(X_{s_{n-1}})V^{\sigma,\lambda_{n-1,j}}_t(h)(s_{n-1},X_{s_{n-1}})
\label{eqmu1}
\end{eqnarray}
From the definition of $\overline h$  on $[s_{n-1},t[\times \R^n$ (equation (\ref{hb1}), it follows that 
$$E_{P^{\sigma, \mu}_{r,y}}(h(X_t)|{\cal B}^r_{s_{n-1}})-\alpha_{s_{n-1}t}(P^{\sigma, \mu}_{r,y}) \leq \overline h(s_{n-1},X_{s_{n-1}})$$
Using the precedings equations and the monotonicity of the conditional expectation we then prove recursively that for all $r \leq s <t$, 
$$E_{P^{\sigma,  \mu}_{r,y}}(h(X_t)|{\cal B}^r_s)-\alpha_{st}(P^{\sigma, \mu}_{r,y}) \leq \overline h(s,X_s)$$
\end{itemize}

\hfill $\square $

 Recall that $(\tilde {\cal Q}_{r,y})_S(\Lambda)$ is the stable set generated by the $\lambda$ continuous bounded $\Lambda$ valued.

\begin{theorem}
 Let $\sigma$ be continuous bounded and $a=\sigma \sigma^*$ strictly elliptic.   Let $0 \leq r \leq s\leq t$. Let $y$ in $\R^n$. Let $\Pi^{r,y}_{s,t}$ and $\alpha_{s,t}$ be defined as in Section \ref{sectionTCC}.
  $$\Pi^{r,y}_{s,t}(Y)=\esssup_{P^{\sigma,  \mu}_{r,y} \in (\tilde {\cal Q}_{r,y})_S(\Lambda)}  (E_{P^{\sigma,  \mu}_{r,y}}(Y|{\cal B}^r_{s})- \alpha_{s,t}({P^{\sigma,  \mu}_{r,y}}))$$
with $\alpha_{s,t}({P^{\sigma, \mu}_{r,y}})=E_{P^{\sigma, \mu}_{r,y}}( \int\limits_{s}^{t} g(u,X_u(\omega),\mu(u,\omega)) d u|{\cal B}^r_{s})$.\\
Assume that $g$ satisfies hypothesis $H_g$.
Let $h$ be a Borelian  map  on $ \R^n$ bounded from below. 
There is a lower semi-continuous  function $\tilde h$ on   $[0,t[ \times \R^n$  such that 
\begin{equation}
\forall s \in [0,t[,\; \forall x \in \R^n,\;\Pi^{s,x}_{s,t}(h(X_t))=\tilde h(s,x)
\label{eqmeas2}
\end{equation}
\begin{equation}
\forall 0 \leq r \leq s <t,\;\;\Pi^{r,y}_{s,t}(h(X_t))=\tilde h (s,X_{s})\; Q^a_{r,y}\;a.s.
\label{eqmeas}
\end{equation}
 In case where $h$ is continuous bounded from below, the map $\tilde h $ extended by  $\tilde h (t,x)=h(x)$ is lower semi-continuous on $[0,t]\times\R^n$. 
\label{thmmar}
\end{theorem}
{\bf Proof}
- Let  $h$ be Borelian  bounded (resp  $h$ continuous bounded). Given $\mu$ in $\tilde{\cal S}^r_t(\Lambda)$ denote $h_{\mu}$ the continuous function on $[r,t[ \times \R^n$ (resp on $[r,t] \times \R^n$) with $h_{\mu}(t,x)=\overline h$ constructed in Theorem \ref{thmMar} satisfying equations (\ref{eqin1}) and (\ref{eqin2}). \\
Notice that for $r \leq s$, every $\nu$ in $\tilde {\cal S}^s_t(\Lambda)$ can be identified with the process $\tilde \nu$ in $\tilde {\cal S}^r_t(\Lambda)$ defined by  $\tilde \nu(u,\omega)=0$ if $u \leq s$ and  $\tilde \nu(u,\omega)=\nu(u,\omega)$ if $s < u$. Let 
\begin{equation}
\tilde h(s,x)=\sup_{\mu \in \tilde {\cal S}^0(\Lambda)} h_{\mu}(s,x)
\label{eqhcont}
\end{equation}
The semi continuity properties for $\tilde h$ follow from the continuity properties for every $h_{\mu}$.
From Theorem \ref{thmMar}, equation (\ref{eqin1}), for given $s$,
\begin{equation}
\tilde h(s,x) \leq \Pi^{s,x}_{s,t}(h(X_t))=\sup_{\nu \in \tilde {\cal S}^s_t(\Lambda)}( E_{P^{\sigma, \nu}_{s,x}}(h(X_t))-\alpha_{st}(P^{\sigma, \nu}_{s,x}))
\label{eqhtilde}
\end{equation}
Furthermore,
\begin{equation}
\sup_{\nu \in \tilde {\cal S}^s_t(\Lambda)} (E_{P^{\sigma, \nu}_{s,x}}(h(X_t))-\alpha_{st}(P^{\sigma, \nu}_{s,x}))
 \leq \sup_{\nu \in \tilde {\cal S}^s_t(\Lambda)}h_{\nu}(s,x) \leq \sup_{\nu \in \tilde {\cal S}^0_t(\Lambda)}h_{\nu}(s,x)=\tilde h(s,x)
\label{eqhtilde1}
\end{equation}
where the first inequality is due to  equation (\ref{eqin2})
 and the second one to the inclusion $\tilde {\cal S}^s_t(\Lambda) \subset \tilde {\cal S}^0_t(\Lambda)$.  This proves equation (\ref{eqmeas2}).
Furthermore from equation (\ref{eqin2}), for all $0 \leq r \leq s$,
\begin{equation}
\Pi^{r,y}_{s,t}(h(X_t))=\sup_{\mu \in \tilde {\cal S}^r_t(\Lambda)}h_{\mu}(s,X_s)
\label{eqhtilde2}
\end{equation}
Notice also that from the inclusions  $\tilde {\cal S}^s_t(\Lambda) \subset 
 \tilde {\cal S}^r_t(\Lambda) \subset \tilde {\cal S}^0_t(\Lambda)$, and the inequalities (\ref{eqhtilde}) and (\ref{eqhtilde1}) it follows that for $0 \leq r \leq s$,
\begin{equation}
\sup_{\mu \in \tilde {\cal S}^r_t(\Lambda)}h_{\mu}(s,x)=\tilde h(s,x)
\label{eqhtilte3}
\end{equation}
Equation (\ref{eqmeas}) follows from (\ref{eqhtilde2})  and (\ref{eqhtilte3}). \\
- For $h$ Borelian bounded from below, or continuous bounded from below, the result follows from the equation 
$$\Pi_{s,t}(h(X_t))=\sup_{n \in\N^*} \Pi_{s,t} (h_n(X_t))$$
with $h_n(x)=h(x)\wedge n$ (cf Corollary \ref{corop2}).
\hfill $\square $


\section{Viscosity solution of the PDE}

\subsection{Continuous selector }
\label{contsec}

Recall the following definition of a continuous selector (Definition 16.57 of \cite{AB}, where a multivalued mapping is called correspondence).

\begin{definition}
A selector from a  multivalued mapping  $\phi$  from  X into Y is a function $s: X \rightarrow Y$ such that $s(x) \in \phi(x)$ for all $x \in X$.  A continuous selector is  a selector  which is continuous.
\label{defsc}
\end{definition}
Recall  also the following definition  from \cite{AB} (Definition 16.2 and Lemma 16.5):

\begin{definition}
A multivalued mapping  $\phi$ from  X into Y is lower hemicontinuous if it satisfies the following equivalent conditions
\begin{itemize}
\item  For every closed subset F of Y, $\phi^u(F)=\{x \in X :  \phi(x) \subset F\}$ is closed
\item  For every open  subset V of Y, $\phi^l(V)=\{x \in X :  \phi(x) \cap V\neq \emptyset\}$ is open
\end{itemize}
\label{deflhc}
\end{definition}
Recall also the following Michel Selection  Theorem (cf \cite{AB} Theorem 16.61)
\begin{theorem}
A lower hemicontinuous mapping from a paracompact space into a Banach space with non empty closed convex values admits a continuous selector.
\label{selcont}
\end{theorem}
Recall also that every metrizable space is paracompact (Theorem 2.86 of \cite{AB}).

\subsection{Comment on the notion of viscosity solution}

\label{seccomvis}

We refer to \cite{GIL} for an exposition of the theory of viscosity solutions for second order partial diffrential equations.
Consider the following PDE on $[0,t[ \times \R^n$
\begin{equation}
H(u,x,v(u,x),\partial_uv(u,x),Dv(u,x),D^2v(u,x))=0
\label{eqPDEH}
\end{equation}
Recall the definition:
\begin{definition} 
\begin{itemize}
\item
An upper semi-continuous function  $v$ is a subsolution in the viscosity sense to (\ref{eqPDEH}) on $[0,t[ \times \R^n$ if for all  $(t_0,x_0)$, $t_0<t$,  and  all  function $\phi$ of class ${\cal C}^{1,2}_b$ such that $(t_0,x_0)$ is a local maximizer of $v-\phi$, and $v(t_0,x_0)=\phi(t_0,x_0)$
\begin{equation}
 H(t_0,x_0,\phi(t_0,x_0),\partial_u \phi(t_0,x_0),D\phi(t_0,x_0),D^2\phi(t_0,x_0))\leq 0
\label{eqsubsol}
\end{equation}
\item A lower semicontinuous function  $v$ is a supersolution in the viscosity sense to (\ref{eqPDEH}) on $[0,t[ \times \R^n$ if for all  $(t_0,x_0)$, $t_0<t$,  and  all  function $\phi$ of class ${\cal C}^{1,2}_b$ such that $(t_0,x_0)$ is a local minimizer of $v-\phi$, and $v(t_0,x_0)=\phi(t_0,x_0)$
\begin{equation}
 H(t_0,x_0,\phi(t_0,x_0),\partial_u \phi(t_0,x_0),D\phi(t_0,x_0),D^2\phi(t_0,x_0))\geq 0
\label{eqsursol}
\end{equation}
\item A continuous function $v$ is a viscosity solution of (\ref{eqPDEH})  on $[0,t[ \times \R^n$ if it is both a subsolution and a supersolution.
\end{itemize}
\end{definition}
In the case where the continuity of $v$ is not known, the semicontinuous envelopes of $v$ are considered:
\begin{equation}
v_*(s,x)=\liminf_{(s',x') \rightarrow (s,x)}v(s',x')
\label{eqliminf}
\end{equation}
\begin{equation}
v^*(s,x)=\limsup_{(s',x') \rightarrow (s,x)}v(s',x')
\label{eqliminf2}
\end{equation}
Assume now that we want to prove that $v$ is a viscosity solution of (\ref{eqPDEH}). We have the following result proving the continuity of $v$ at all points as soon as $v$ is continuous at terminal points, when a comparison principle is satisfied.
\begin{proposition}
Let $v: [0,t] \times \R^n \rightarrow \R$. Let $g(x)=v(t,x)$. Assume that $v$ is  continuous at $(t,x)$ for all $x \in \R^n$.\\
Assume that  $v_*$ is a supersolution of (\ref{eqPDEH}) in the viscosity sense on $[0,t[ \times \R^n$. Assume that $v^*$ is a subsolution of (\ref{eqPDEH}) in the viscosity sense on $[0,t[ \times \R^n$. 
Assume that  a comparison principle is satisfied. \\
Then $v=v^*=v_*$, $v$ is continuous on $[0,t] \times \R^n$ and is the unique solution of viscosity of 
(\ref{eqPDEH}) on $[0,t[ \times \R^n$ with terminal condition $g$.
\label{propvis}
\end{proposition}
 {\bf Proof} 
The function $v$ being continuous at $(t,x)$ for all $x$ in $\R^n$, it follows that $v^*(t,x)=v_*(t,x)=v(t,x)=g(x)$ for all $x$.
The function $v_*$ being a supersolution, $v^*$ a subsolution, it follows then from the comparison principle  that $v^* \leq v_*$.\\
On the other hand it  follows from the definition of $v_*$ and  $v^*$ that $v_* \leq v \leq  v^*$.
This proves that $v_*=v=v^*$. Thus  $v$ is continuous and  is the unique viscosity solution of  (\ref{eqPDEH})  on $[0,t[ \times \R^n$ with terminal condition $g$. \hfill $\square $  

\begin{remark} 
 If $v$ is not continuous at $(t,x)$ for some $x$ $v_*(t,x)< v^*(t,x)$. The comparaison principle does not allow to conclude.
\end{remark}
\subsection{Viscosity supersolution }
\label{secsuper}
We introduce the following hypothesis:
\begin {definition}
Hypothesis $H_{\Lambda}$:\\
The multivalued Borel mapping $\Lambda$ satisfies hypothesis  $ H_{\Lambda}$ if
  $\Lambda$ is a multivalued  Borel mapping convex and closed valued such that for all $K $ large enough, $\Lambda_K$ is lower hemicontinuous, where $\Lambda_K(t,x)=\{y \in \Lambda(t,x),\;||y|| \leq K\}$.
\label{defHl}
\end{definition}

Assume that $h$ is continuous   bounded from below and that $g$ satisfies hypothesis $H_g$. Let $t>0$. In Theorem \ref{thmmar}, we have proved the existence of  a lower semi continuous function $v=\tilde h$ on $[0,t] \times \R^n$  such that 
\begin{equation}
v(s,x)=\Pi^{s,x}_{s,t}(h(X_t))
\label{eqv1}
\end{equation}
 and 
\begin{equation}
\Pi^{r,y}_{s,t}(h(X_t))=v(s,X_s)\;\;  \forall r \leq s\leq t\;\;\; Q^a_{r,y}\;a.s.
\label{eqv2}
\end{equation}

We want to prove that $v$ is a viscosity supersolution of a second order PDE.

 Let $(t_0,x_0) \in [0,t[ \times \R^n$.
 Let $\phi \in {\cal C}^{1,2}_b([0,t] \times \R^n )$  such that 
$$0=v(t_0,x_0)-\phi(t_0,x_0)=\min(v(t,x)-\phi(t,x))$$

\begin{lemma} Assume that $\sigma^{-1}(u,X_u)$ is bounded on $[0,T] \times \Omega$. For all bounded $\mu$ such that $P^{\sigma,  \mu}_{r,y}$ belongs to $\tilde {\cal Q}{r,y}(\Lambda)$, there is a Brownian motion $W^{\mu}(u)$ with respect to $P^{\sigma,  \mu}_{r,y}$ such that 
\begin{equation}
dX_u(\omega)=\sigma(u,X_u)(\omega)\mu(u,\omega)du + \sigma(u,X_u)(\omega) d W^{\mu}(u)(\omega)\;\;\;P^{\sigma, \mu}_{r,y}\;a.s.
\label{eqBrown}
\end{equation}
\label{lemmaBrown}
\end{lemma}
{\bf Proof} From Proposition \ref{uniqmar},  $P^{\sigma, \mu}_{r,y}$ is  solution to the martingale problem
\begin{equation}
Y_t=exp \{(\theta)^*(X(t)-y)-\int_r^t \theta^*\sigma(u, X_u(\omega)) \mu(u,\omega)  du -\frac{1}{2} \int_r^t \theta^* a(u,X_u(\omega)) \theta du \} 
\label{eqmarmu0}
\end{equation}
starting from $y$ at time $r$.\\
Let $\theta \in \R^n$. Since the functions $\sigma$ and $\sigma^{-1}$ are bounded, and the process $\mu$ is bounded and non anticipating, one can apply  Theorem 3.2 of \cite{SV1}  with $$\theta(u,\omega)=(\sigma^{-1})^*(u,\omega)\theta$$
and $$\xi(s)(\omega)=X_s(\omega)-y-\int_r^s\sigma(u,X_u(\omega))\mu(u,\omega)du$$
It follows that 
$$exp(\int_r^t \theta^*(\sigma^{-1})(u,\omega)[dX_u-\sigma(u,X_u(\omega))\mu(u,\omega)du])-\frac{1}{2}|\theta|^2(t-s)$$
is a  $P^{\sigma, \mu}_{r,y}$  martingale. The proof follows then the proof of Theorem 3.3 of \cite{SV1}. \hfill $\square $ 

\begin{theorem}Let $\sigma$ be continuous bounded and $a=\sigma \sigma^*$ strictly elliptic.   Assume that $\Lambda$  satisfies  hypothesis $H_{\Lambda}$. Let $h$ be  continuous   bounded from below. Assume that  the restriction of $g$  to $\{(u,x,y),\; y \in \Lambda(u,x)\}$ is upper semi-continuous   and that $g$ satisfies hypothesis $H_g$.  Let $v$ be a lower semi continuous function such that equations (\ref{eqv1}) and (\ref{eqv2}) are satisfied. Then $v$ is a viscosity supersolution of 
\begin{equation}
\left \{\parbox{12cm}{
\begin{eqnarray} 
-\partial_u v(u,x)-{\cal L} v(u,x)-f(u, x,\sigma^*(u,x)Dv(u,x))&=&0 \nonumber \label{eqedp00}\\
v(t,x)& = & h(x)\nonumber 
\end{eqnarray}
}\right.
\end{equation}
at each point $(t_0,x_0)$ such that  $f(t_0,x_0, \sigma(t_0,x_0)^*D\phi(t_0,x_0))<\infty$.
Here $f(u,x,z)=\sup_{\lambda \in \Lambda(u,x)}(z^*\lambda-g(u,x,\lambda))$,  and 
$ {\cal L}v(u,x) = \frac{1}{2} Tr(a(t,x)D^2(v)(t,x))$.

\label{thmsuper}
\end{theorem}
{\bf Proof}
\begin{itemize}
\item  Step 1: Time consistency\\
From the time consistency property for $\Pi^{t_0,x_0}_{u,t}$, it follows that for all stopping time $0<\delta < t-t_0$ taking a finite number of real values,
\begin{eqnarray}
\Pi^{t_0,x_0}_{t_0,t}(h(X_t))=\Pi^{t_0,x_0}_{t_0,t_0+\delta}(\Pi^{t_0,x_0}_{t_0 +\delta,t}(h(X_t)))\;\nonumber \\
=\sup_{P^{\sigma, \mu}_{t_0,x_0} \in ( \tilde {\cal Q}_{t_0,x_0})_S(\Lambda)}  (E_{P^{\sigma, \mu}_{t_0,x_0}}(\Pi^{t_0,x_0}_{t_0 +\delta,t}(h(X_t)))- \alpha_{t_0,t_0+\delta}({P^{\sigma, \mu}_{t_0,x_0}}))
\label{eqtc}
\end{eqnarray}
Recall that $v$ is lower semi continuous, $\phi$ is continuous, and  $v(u,x) \geq \phi(u,x)$  for all $u \in [0,t]$ and $x \in \R^n$. From equation (\ref{eqv2}) it then follows that  for every stopping time $\delta$ taking a finite number of real  values, 
\begin{equation}
\Pi^{t_0,x_0}_{t_0+\delta,t}(h(X_t))=v(t_0+\delta, X_{t_0+\delta}) \geq \phi(t_0+\delta, X_{t_0+\delta})
\label{ineq02}
\end{equation}
\item Step 2: Ito's formula in case $\sigma^{-1}$ is bounded.\\
 From Lemma \ref{lemmaBrown} applied with $\mu(u,\omega)=\lambda(u,X_u(\omega))$, there is a   Brownian motion  $W^{\lambda}$  with respect to  the probability measure $Q^{\sigma, \lambda}_{t_0,x_0}$ such that 
\begin{equation}
dX_u=(\sigma\lambda)(u,X_u)du + \sigma(u,X_u) dW^{\lambda}_u\;\;\
\label{eqBrown2}
\end{equation}
We can now apply Ito's formula to  $\phi(t,X_t)$.
\begin{eqnarray}
\phi(t_0+\delta, X_{t_0+\delta})=\phi(t_0,x_{0})\nonumber\\+ \int_{t_0}^{t_0+\delta}[\phi_u(u, X_u)+ \frac{1}{2}Trace(D^2\phi(u,X_u)(a)(u,X_u))]du \nonumber\\
+ \int_{t_0}^{t_0+\delta}(D\phi)^*(u, X_u)\sigma(u,X_u)dW_u^{\lambda}
+ \int_{t_0}^{t_0+\delta}(D\phi)^*(u, X_u)\sigma(u,X_u)\lambda(u,X_u)du
\label{ineq3}
\end{eqnarray}
It follows from equations  (\ref{eqtc}),  (\ref{ineq02}) and (\ref{ineq3}), that for every $Q^{\sigma, \lambda}_{t_0,x_0}$ ,
\begin{eqnarray}
\Pi^{t_0,x_0}_{t_0,t}(h(X_t))-\phi(t_0,x_{0})
\geq \nonumber \\
E_{Q^{\sigma, \lambda}_{t_0,x_0}}
[\Pi^{t_0,x_0}_{t_0+\delta,t}(h(X_t))
-\int_{t_0}^{t_0+\delta}g(u, X_u,\lambda(u,X_u)du] -\phi(t_0,x_{0})\geq 
\nonumber\\
E_{Q^{\sigma, \lambda}_{t_0,x_0}}[\int_{t_0}^{t_0+\delta}(\phi_u(u, X_u)+ \frac{1}{2} Trace(D^2\phi(u,X_u)a(u,X_u))du]\nonumber\\
+E_{Q^{\sigma, \lambda}_{t_0,x_0}}[ \int_{t_0}^{t_0+\delta}(D\phi^*(u, X_u)\sigma(u,X_u)\lambda(u,X_u)-g(u, X_u,\lambda(u,X_u))
du]\nonumber\\
\end{eqnarray}
Thus
\begin{eqnarray}
0 \geq 
{ E_{Q^{\sigma, \lambda}_{t_0,x_0}}[\int_{t_0}^{t_0+\delta}(\phi_u(u, X_u)+ \frac{1}{2}Trace(D^2\phi(u,X_u)a(u,X_u))du}\nonumber\\
+{ E_{Q^{\sigma, \lambda}_{t_0,x_0}}[ \int_{t_0}^{t_0+\delta}(D\phi^*(u, X_u)\sigma(u,X_u)\lambda(u,X_u)-g
(u, X_u,\lambda(u,X_u))
du}\nonumber\\
\label{eqvs}
\end{eqnarray}
\item Step 3 General case \\
We consider  $\sigma _n$ such that 
  $(\sigma_n)^{-1}$ is bounded and such that the sequence $\sigma_n$ is uniformly bounded and converges to $\sigma$ uniformly on compact spaces.  From Theorem 9.2 of \cite{SV2} applied for given $\lambda$, $t_0$ and $x_0$, $Q^{\sigma_n, \lambda}_{t_0,x_0}$ converges weakly to $Q^{\sigma, \lambda}_{t_0,x_0}$  as $n \rightarrow \infty$.  The function $$\int_{t_0}^{t_0+\delta}(\phi_u(u, X_u)+ \frac{1}{2}Trace((D^2\phi)\sigma^2)(u,X_u))+((D\phi)^*\sigma\lambda)(u,X_u))du$$ being a continuous bounded function of $\omega$, we can pass to the limit when $n \rightarrow \infty$.
$g$ is not bounded but $\lambda$ is bounded. Using the polynomial growth hypothesis for $g$ and the estimates of the moments of $X_t$ (Proposition  \ref{propKryl})  we can proceed as in the proof of Proposition \ref{propcont3}.
Let $g_k=\sup(\inf(g,k),-k)$. 
There is $\tilde k>0$ such that $\forall k \geq \tilde k$, for all $\tilde \sigma$= $\sigma_n$ or $\sigma$,
\begin{equation}
E_{Q^{\tilde\sigma, \lambda}_{t_0,x_0}}(\int_{t_0}^{t_0+\delta}|g(u,X_u, \lambda(u,X_u))-g_k(u,X_u, \lambda(u,X_u))|)du \leq \epsilon
\label{eqinte}
\end{equation}
Thus we can also pass to the limit as $\sigma_n$ tends to $\sigma$ in equation (\ref{eqvs}) even for the $''g''$ term.
This proves that inequality  (\ref{eqvs}) is satisfied for all $\sigma$ bounded such that $\sigma \sigma^*$ is strictly elliptic. 
\item last step: viscosity supersolution\\
By hypothesis $f(t_0,x_0, \sigma(t_0,x_0)^*D\phi(t_0,x_0))< \infty$.
Thus for all  $\epsilon>0$, there is $\lambda_0 \in \Lambda(t_0,x_0)$ such that 
\begin{equation}
((D\phi)^*\sigma)(t_0,x_0)\lambda_0-g(t_0,x_0,\lambda_0)> f(t_0,x_0, \sigma(t_0,x_0)^*D\phi(t_0,x_0)) - \epsilon
\label{eqlambda0}
\end{equation}
From Hypothesis ($H_{\Lambda}$) for $K$ large enough $\Lambda_K$ is lower hemicontinuous. Choose such a $K>||\lambda_0||$. Let $C$ be the multivalued Borel mapping defined by 
\begin{eqnarray}
C(u,x)&=&\Lambda_K(u,x),\;\; \forall (u,x) \neq (t_0,x_0)\nonumber\\
C(t_0,x_0)&=&\{\lambda_0\}\nonumber
\end{eqnarray}
$C$ is also lower hemicontinuous.
 From Theorem \ref{selcont},   there is thus a continuous bounded selector   $\lambda(u,x)$  of $\Lambda$ such that for every $(u,x)$, $\lambda(u,x) \in \Lambda_K(u,x)$ and $\lambda(t_0,x_0)=\lambda_0$. \\
From the upper semi-continuity of $g$ on $\{(u,x,y),\; y \in \Lambda(u,x)\}$, and the continuity  of the maps $\lambda$,  $X$ $\Phi_u$,$D^2\Phi$, and $\sigma$, for all $\epsilon>0$, there is $\eta>0$,  such that for $t_0 \leq u \leq t \leq t_0+\eta$ and $||x_0-x||<\eta$, 
\begin{equation}
g(t,x,\lambda(t,x))-g(t_0,x_0,\lambda_0)< \epsilon
\label{eqg0}
\end{equation} 
\begin{eqnarray}
|\phi_u(u, x)+ \frac{1}{2}Trace(D^2\phi(u,x)a(u,x))+(D\phi^*(u, x)\sigma(u,x)\lambda(u,x)-\nonumber\\\phi_u
(t_0
, x_0)+ \frac{1}{2}Trace(D^2\phi(t_0,x_0)
(a)(t_0,x_0))+D\phi^*(t_0,x_0)
\sigma(t_0,x_0)\lambda_0| \leq \epsilon
\label{eqsubv}
\end{eqnarray}
From Proposition \ref{propKryl}, applied to the probability measure ${Q^{\sigma, \lambda}_{t_0,x_0}}$, there is $0<\alpha < \eta$  such that 
$${Q^{\sigma, \lambda}_{t_0,x_0}}(A)<\epsilon\;\; with \;\;A=\{\omega \;|\; \sup_{t_0 \leq u \leq t_0+ \alpha}||X_u-x_0||>\eta\;\}$$

Let $\delta=\alpha 1_{A^c}$. $\delta$ is a stopping time taking only two values.
For $t_0 \leq u \leq t_0+\delta$, $g(u,X_u(\omega),\lambda(u,X_u(\omega)))<g(t_0,x_0,\lambda_0)+\epsilon$.
Divide the inequality (\ref{eqvs}) by $\alpha$. The left hand side is equal to $0$. From equations (\ref{eqg0}),  (\ref{eqsubv}) and  (\ref{eqlambda0}), the  right hand side is greater or equal to 
$[\phi_u
(t_0
, x_0)+ \frac{1}{2}Trace(D^2\phi(t_0,x_0)
(a)(t_0,x_0))+ f(t_0,x_0, \sigma(t_0,x_0)^*D\phi(t_0,x_0))-3\epsilon](1-\epsilon)$.
This proves that 
$v$ is a viscosity supersolution of (\ref{eqedp00}).
\end{itemize}
\hfill $\square $ 
\subsection{Viscosity subsolution}
As in the previous subsection we assume that $h$ is continuous and bounded from below.
As in the previous Section $v=\tilde h$ is the lower semicontinuous function   on $[0,t] \times \R^n$  such that $v(t,x)=h(x)$ and $\Pi^{s,x}_{s,t}(h(X_t))=v(s,x)$  and 
\begin{equation}
\Pi^{r,y}_{s,t}(h(X_t))=v(s,X_s)\;\;\forall \;\;r \leq s\leq t
\label{eqFellsc}
\end{equation}
 (cf Theorem \ref{thmmar}).

Denote $v^*$ the upper semi continuous envelope of $v$, $$v^*(s,x)=\limsup_{(s',x')\rightarrow (s,x)}    v(s',x')$$
Let $\phi \in {\cal C}^{1,2}_{b}$ such that 
\begin{equation}
0=v^*(t_0,x_0)-\phi(t_0,x_0)=\sup_{(s,x)}(v^*(s,x)-\phi(s,x))
\label{equsc}
\end{equation}

We introduce now a growth hypothesis on $\Lambda$.
\begin{definition}

\begin{itemize}
\item The multivaled Borel mapping $\Lambda$ has linear growth if 
 there is a Borelian  map $\phi_{\Lambda}$  such that 
\begin{equation}
\forall y \in \Lambda(t,x),\; ||y|| \leq \phi_{\Lambda}(t,x)
\label{eqphi}
\end{equation}
and such that  $\phi_{\Lambda}(s,x) \leq K(1+||x||)$.
\item The multivaled Borel mapping $\Lambda$   satisfies hypothesis $\overline H_{\Lambda}$ if it satisfies hypothesis  $H_{\Lambda}$ (cf Definition \ref{defHl}) and if it has linear growth.
\end{itemize}
\end{definition}

\begin{lemma}
Assume that   the multivaled Borel mapping $\Lambda$ has linear growth. For all $q\geq 1$, $A$, $C$,  and $t>0$, there are  constants $K_1$ such that for all $y$ such that $||y|| \leq C$ and $\sigma$ continuous such that  $||\sigma|| \leq A$, and all bounded $\mu$ such that $P^{\sigma, \mu}_{r,y}$ belongs to $(\tilde{\cal Q}_{r,y})_S(\Lambda)$
\begin{equation}
E_{P^{\sigma, \mu}_{r,y}}(\sup_{s \leq u \leq t }(||X_t-X_{u}||^{2q}) \leq K_1 (t-s)^q
\label{eqKry01}
\end{equation}
\label{lemKrylb}
\end{lemma}
{\bf Proof}
\begin{itemize}
\item
Assume that $\sigma^{-1}(u,X_u)$ is bounded on $[0,T] \times \Omega$. From Lemma \ref{lemmaBrown}, for all bounded $\mu$ such that $P^{\sigma, \mu}_{r,y}$ belongs to $\tilde {\cal Q}_{r,y}(\Lambda)$, there is a Brownian motion $W^{\mu}_u$ with respect to $P^{\sigma, \mu}_{r,y}$ such that 
\begin{equation}
dX_u(\omega)=\sigma(u,X_u)(\omega)\mu(u,\omega)du + \sigma(u,X_u)(\omega) d W^{\mu}_u(\omega)\;\;\;P^{\sigma, \mu}_{r,y}\;a.s.
\label{eqBrown1}
\end{equation}
Let $b(s,X_s(\omega),\omega)=\sigma(s,X_s(\omega))\mu(s,\omega)$.  $\sigma$ being bounded by $A$, and $\mu$ being $\tilde \Lambda$ valued, it follows  from hypothesis $\overline H_{\Lambda}$ that  $||b(s,X_s(\omega),\omega)|| \leq AK(1+||X_s(\omega)||)$.
 It follows then  from  \cite{Kr}, II 5 Corollary 10, that there exists  $K_1>0$  depending only on $q$ $A$ $B$ $K$ and $t$ such that  equation (\ref{eqKry01}) 
is satisfied.
\item For general $\sigma$ consider as in  Theorem \ref{thmsuper} a sequence $\sigma _n$ such that $(\sigma_n)^{-1}$ is bounded and such that the sequence $\sigma_n$ is uniformly bounded and converges to $\sigma$ uniformly on compact spaces.
We then conclude as in the proof of Proposition \ref{propKryl}.
\end{itemize}
\hfill $\square $

\begin{theorem}  Let $\sigma$ be continuous bounded and $a=\sigma \sigma^*$ strictly elliptic.   Assume that $\Lambda$  has linear growth. Assume that $h$ is continuous bounded from below. Assume that $g$  satisfies hypothesis $H_g$. Assume that $f$ is upper semi continuous at $(t_0,x_0,\sigma^*(t_0,x_0)D\phi(t_0,x_0))$.  
Let $v=\tilde h$ be the lower semi continuous function  as in Theorem \ref{thmmar} . Then $v^*$ is a viscosity subsolution of 
\begin{equation}
-\partial_u v(u,x)- {\cal L} v(u,x)-f(t, x,\sigma^*(u,x)Dv(u,x))=0 
 \label{eqedp01}
\end{equation}
with $f(u,x,z)=\sup_{\lambda \in \Lambda(u,x)}(z^*\lambda-g(u,x,\lambda))$
\label{thmsub}
\end{theorem}

{\bf Proof}
The function $f$ being upper semi continuous at  $(t_0,x_0,\sigma^*(t_0,x_0)D\phi(t_0,x_0))$, and $X$, $\Phi_u$, $D^2\Phi$, $D\phi$  and $\sigma$ being continuous, for all $n \in \N^*$, there is $\eta_n>0$, $\eta_n<1$, such that for $t_0 \leq u \leq t \leq t_0+\eta_n$ and $||x_0-x||<\eta_n$, 
 \begin{equation}
f(t,x,\sigma^*(t,x)D\phi(t,x)) \leq   f(t_0,x_0,\sigma^*(t_0,x_0)D\phi(t_0,x_0))+ \frac{1}{n}
\label{eqviscsub}
\end{equation}
\begin{eqnarray}
|\phi_u(u, x)+ \frac{1}{2}Trace(D^2\phi(u,x)a(u,x))+(D\phi^*(u, x)\sigma(u,x)\lambda(u,x)-\nonumber\\\phi_u
(t_0
, x_0)+ \frac{1}{2}Trace(D^2\phi(t_0,x_0)
a(t_0,x_0))+D\phi^*(t_0,x_0)
\sigma(t_0,x_0)\lambda_0| \leq \frac{1}{n}
\label{eqsupv}
\end{eqnarray}
From Lemma \ref{lemKrylb}, there is $\gamma_n>0$ such that for all $|u-t_0|\leq 1$ and $||y-x_0|| \leq 1$ for all $\mu$ $\tilde \Lambda$ valued such that 
$P^{\sigma, \mu}_{u,y}$ belongs to $\tilde {\cal Q}_{u,y}(\Lambda)$,
$${P^{\sigma, \mu}_{u,y}}(A_n)<\frac{1}{n}\;\; \text{with} \;\;A_n=\{\omega \;|\; \sup_{u \leq u' \leq u+ \gamma_n}||X_{u'}-y||>\frac{\eta_n}{2}\;\}$$
Let $\delta_n=\gamma_n 1_{A_n^c}$. $\delta_n$ is a stopping time taking only two values.\\
For all $n>0$ choose $(t_n,x_n)$ such that  $|t_n-t_0|<\frac{\eta_n}{2}$, $||x_n-x_0||<\frac{\eta_n}{2}$ and $\phi(t_n,x_n)-  v(t_n,x_n) \leq \frac{\gamma_n}{n}$. 
\\
From equation (\ref{eqFellsc}) $v(t_n,x_n)=\Pi^{t_n,x_n}_{t_n,t}(h(X_t))$. On the other hand,
\begin{equation}
\Pi^{t_n,x_n}_{t_n,t}(h(X_t))
={\sup}_{P^{\sigma, \mu}_{t_n,x_n} \in (\tilde {\cal Q}_{t_n,x_n})_S(\Lambda)}(E_{P^{\sigma, \mu}_{t_n,x_n}}(h(X_t))-\alpha_{t_n,t}(P^{\sigma, \mu}_{t_n,x_n}))
\label{eqmun}
\end{equation}
Thus for all $n>0$, there is a bounded process $\mu_n$, $t_n$-non-anticipating,  $P^{\sigma, \mu_n}_{t_n,x_n} \in (\tilde{\cal Q}_{t_n,x_n})_S(\Lambda)$ 
such that 
\begin{equation}
v(t_n,x_n) \leq E_{P^{\sigma, \mu_n}_{t_n,x_n}}(h(X_t))-\alpha_{t_n,t}(P^{\sigma, \mu_n}_{t_n,x_n})+\alpha_n
\label{eq24}
\end{equation}
From the cocycle condition for the penalty associated to the probability measure $P^{\sigma,\mu_n}_{t_n,x_n}$ and the definition of $\Pi^{t_n,x_n}_{t_n+\delta_n,t}(h(X_t))$, it follows from (\ref{eq24}) that 
\begin{equation}
v(t_n,x_n) \leq E_{P^{\sigma, \mu_n}_{t_n,x_n}}[\Pi^{t_n,x_n}_{t_n+\delta_n,t}(h(X_t))-\int_{t_n}^{t_n+\delta_n}g (u,
X_u,\mu_n(u,\omega)du]+\alpha_n
\label{eq25}
\end{equation}
From equation (\ref{eqFellsc}) 
\begin{equation}
\Pi^{t_n,x_n}_{t_n+\delta_n,t}(h(X_t))=v(t_n+\delta_n,X_{t_n+\delta_n})
\label{eqpi}
\end{equation}
Furthermore, $v \leq v^* \leq \phi$, $v$ being lower semi continuous, $v^*$ upper semi continuous  and $\phi$  continuous. It follows that 
\begin{equation}
v(t_n+\delta_n,X_{t_n+\delta_n}) \leq \phi(t_n+\delta_n,X_{t_n+\delta_n})
\label{eqMelies}
\end{equation}
 From equations (\ref{eq25}) (\ref{eqpi}) and (\ref{eqMelies}), it  follows that 
$$v(t_n,x_n) \leq E_{P^{\sigma, \mu_n}_{t_n,x_n}}[\phi(t_n+\delta_n,X_{t_n+\delta_n})-\int_{t_n}^{t_n+\delta_n}g (u,X_u,\mu_n(u,\omega)du]+\alpha_n$$
As in the proof of Theorem \ref{thmsuper} we consider  $\sigma_j$ continuous, $\sigma_j^{-1}$ bounded  such that the sequence  $\sigma_j$ is uniformly boundedby $A$  and uniformly convergent to $\sigma$ on compact spaces. Proceeding as in the proof of Proposition  \ref{propcont3}, even if $g$ is not bounded,  using the weak convergence of  $P^{\sigma_j,\mu_n}_{t_0,x_0}$ to $P^{\sigma, \mu_n}_{t_0,x_0}$  when  $j \rightarrow \infty$, we  get that
for given $n$,  there is $j(n)$ such that 
\begin{equation}
v(t_n,x_n) \leq E_{P^{\sigma_{j(n)}, \mu_n}_{t_n,x_n}}[\phi(t_n+\delta_n,X_{t_n+\delta_n})-\int_{t_n}^{t_n+\delta_n}g (u,X_u,\mu_n(u,\omega)du]+2\alpha_n
\label{eqMelies2}
\end{equation}
As in the the proof of Theorem \ref{thmsuper}, we apply Ito's formula to $\phi(t_n+\delta_n,X_{t_n+\delta_n})$ using the ${P^{\sigma_{j(n)}, \mu_n}_{t_n,x_n}}$-Brownian motion $W^{n}$. 
Using the inequality
\begin{eqnarray}
D\phi^*(u,X_u(\omega))\sigma_{j(n)}(u,X_u(\omega))\mu_n(u,\omega)-g(u,X_u,\mu_n(u,\omega)) \leq \nonumber\\ f(u,X_u(\omega),\sigma_{j(n)}^*(u,X_u(\omega))D\phi(u,X_u(\omega)) \label{eqmu0}
\end{eqnarray}
 Letting $n$ tend to $\infty$, it follows from the definition of $\delta_n$, the uniform convergence of $\sigma_n$ to $\sigma$ on compact spaces and upper semicontinuity of $f$ at $(t_0,x_0,\sigma^*(t_0,x_0)D\phi(t_0,x_0))$ that 
$$0 \leq \phi_u(t_0, x_0)+  {\cal L} \phi(t_0,x_0) + f(t_0,x_0,\sigma(t_0,x_0)D\phi(t_0,x_0))$$ i.e.
$v^*$ is a viscosity subsolution of (\ref{eqedp01}).\\

\hfill $\square $

\subsection{Viscosity solution and uniqueness}
\label{subsecU}
The following Theorem results from Theorems \ref{thmsuper} and \ref{thmsub}
\begin{theorem}
Let $\sigma$ be continuous bounded and $a=\sigma \sigma^*$ strictly elliptic. 
Assume that $g$ satisfies hypothesis $ H_g$. Assume that $\Lambda$  satisfies hypothesis $\overline H_{\Lambda}$. Let $h$ be continuous bounded from below.   Assume that $f$ is upper semi continuous and that the restriction of $g$ to  $\{(u,x,y),\;y \in \Lambda(u,x)\}$ is upper semi continuous.  Let $v=\tilde h$ be the lower semi continuous function on $[0,t] \times \R^n$ as in Theorem \ref{thmmar}. Assume that the function $v$ is continuous. Then $v$   is a viscosity solution on $[0,t[ \times \R^n$ of 
\begin{equation}
\left \{\parbox{12cm}{
\begin{eqnarray} 
-\partial_u v(u,x)-{\cal L} v(u,x)-f(u, x,\sigma^*(u,x)Dv(u,x))&=&0 \nonumber \label{eqedp02}\\
v(t,x)& = & h(x)\nonumber 
\end{eqnarray}
}\right.
\end{equation}
where  $f(u,x,z)=\sup_{\lambda \in \Lambda(u,x)}(z^*\lambda-g(u,x,\lambda)$
\label{thmuniqvis}
\end{theorem}

We now give two sufficient conditions for the continuity of the function $v$. The first one is the existence of an optimal probability measure. The following lemma is a direct application of proposition \ref{propcontfell}.

\begin{lemma}
Assume that there is a probability measure $Q^{ \sigma,\lambda}_{0,y}$  in ${\tilde{\cal Q}}_{0,y}(\Lambda)$  such that for all $s \in [0,t[$, 
\begin{equation}
\Pi^{0,y}_{s,t}(h(X_t)=E_{Q^{\sigma, \lambda}_{0,y}}(f(X_t)|{\cal B}_{s})- \alpha_{s,t}(Q^{\sigma, \lambda}_{0,y})
\label{max}
\end{equation}
Then the function $v=\tilde h$  as in Theorem \ref{thmmar}  is continuous.
\label{lemmaatt}
\end{lemma}

The following lemma gives a sufficient condition for the continuity in cases where an optimal control does not exist.
                                                                                                                                                                
\begin{lemma}
Let $\sigma$ be continuous bounded and $a=\sigma \sigma^*$ strictly elliptic. Assume that $g$ satisfies hypothesis $H_g$. Assume that $\Lambda$ satisfies hypothesis $\overline H_{\Lambda}$.  Let $h$ be bounded from below and $\alpha$ H\"older-continuous for some $\alpha>0$.  Let $v=\tilde h$ be the lower semi continuous function on $[0,t] \times \R^n$ as in Theorem \ref{thmmar}. Assume that the PDE (\ref{eqedp02}) satisfies the comparison principle for functions bounded on compact spaces. Then the function $v$ is continuous. It is the unique viscosity solution of  (\ref{eqedp02}). 
\label{lemmacontt}
\end{lemma}

{\bf Proof} 
\begin{itemize}
\item
We prove first  that $v$ is  continuous at $(t,x)$ for all $x$ and  that $v$ is bounded on compact sets.
Recall that for all $x$, $v(t,x)=h(x)$. Let $s<t$.
\begin{equation} 
|v(s,x)-h(x)| \leq \sup_{\mu} |E_{P^{\sigma, \mu}_{s,x}}(h(X_t)-h(x))|+ \sup_{\mu} |\alpha_{st}(P^{\sigma, \mu}_{s,x})|
\label{eqcont1}
\end{equation}
Let  $m$ be the growth exponent of $g$ and $\alpha>0$ the exponent of H\"older-continuity of $h$,  one can choose $p>1$ such that $\alpha p>2$ and $mp>2$. 
From  H\"older inequality it follows  that there is $\tilde C$ such that 
\begin{eqnarray} 
|v(s,x)-h(x)| \leq \tilde C \sup_{\mu}([E_{P^{\sigma, \mu}_{s,x}}|X_t-x|^{\alpha p}]^{\frac{1}{p}}+(t-s)\nonumber\\
+[E_{P^{\sigma, \mu}_{s,x}}(\int_s^T||X_u||^{mp}du)]^{\frac{1}{p}})
\label{eqcont3}
\end{eqnarray}
It follows from Lemma \ref{lemKrylb} that 
$|v(s,x)-h(x)|$ tends to $0$  uniformly in $x$ when $s$ tends to $t$.
Since $h$ is continuous and $h(x)=v(t,x)$, it follows that $v$ is  continuous at $(t,x)$ for all $x$ and also  that $v$ is bounded on compact sets.  
\item  It follows from the definition of $v^*$ that $v^*(t,x)=v(t,x)=h(x)$ for all $x$.  From Theorem  
\ref{thmsuper} $v$ is a viscosity supersolution of (\ref{eqedp02}). From Theorem \ref{thmsub}, the function $v^*$ is  a viscosity subsolution of (\ref{eqedp02}).  The function $v$, and thus   $v^*$, are bounded on compact spaces. It follows from the comparison principle that $v^* \leq v$ on $[0,T] \times \R^n$. The converse inequality follows from the definition of $v^*$, thus $v=v^*$ is continuous and is the unique viscosity solution of  (\ref{eqedp02}).
\hfill $\square$
\end{itemize}
For comparison results for non linear second order PDE we refer to \cite{GIL,FSo,Ko,Da}.\\
We give now an application to a stochastic volatility model.


\section{Application to a Stochastic volatility Model}
\label{stocvol}
In this section we consider a general framework for a stochastic volatility model in mathematical finance. The volatility process for the price process is assumed to be itself a stochastic process. The financial market is thus incomplete and we assume that it satisfies the no arbitrage hypothesis. Therefore there exists an equivalent martingale measure for the discounted price process $S_t$.
A general framework for a stochastic volatility model written under an equivalent martingale  measure is thus: 

\[
\left \{ \parbox{12cm}{
\begin{eqnarray} 
dS_t& = & \sigma(t,S_t,Y_t)(\sqrt{1-\rho(t,S_t,Y_t)^2}dW^1_t+\rho (t,S_t,Y_t)dW^2_t) \nonumber\\
dY_t & = & \alpha(t,S_t,Y_t)dt+\gamma(t,S_t,Y_t)dW^2_t\nonumber
\end{eqnarray}
}\right.
\]
where $W^1$ and $W^2$ are two independent n dimensional Brownian motions.
This framework includes the framework considered by Hull and White \cite{HW}.
Assume that $\sigma$ and $\gamma$ are  continuous bounded symmetric definite positive matrices  for all $(t,s,y)$, and that   $\rho$ is  continuous with  values in $]-1,1[$. 
Denote $\Sigma$ the matrix associated to the above system, and $a=\Sigma \Sigma^*$. 
\begin{eqnarray}
& a(t,S_t,Y_t)=\nonumber\\
& \left(
\begin{array}{ccc}
\sigma^2(t,S_t,Y_t) & (\rho\sigma  \gamma)(t,S_t,Y_t)) \\
( \rho )\gamma \sigma(t,S_t,Y_t)) & (\gamma^2(t,S_t,Y_t)
\end{array}
\right)
\end{eqnarray}
One can  easily   verify that $a$ is continuous bounded strictly elliptic.
Define now  the multivalued Borel mapping 
$\Lambda$ from $\R_+ \times \R^{2n}$ into $\R^{2n}$.
\begin{equation}
\Lambda(t,s,y)=\{(\alpha,\nu),\;|\;\alpha \sqrt{(1-\rho(t,s,y)^2)}+\nu \rho(t,s,y)=0\}
\label{eqdeflambda}
\end{equation}
It follows easily from the definition that for all $(t,s,y)$, $\Lambda(t,s,y)$ is a $n$ dimensional vector subspace of $\R^{2n}$. It is thus convex and closed.  A  probability measure $P^{\Sigma,\mu}_{S_0,Y_O}$  is an  equivalent martingale measures for $S_t$ if and only if for all $t$ and $\omega$, $\mu_t(\omega)$ belongs to  $\Lambda(t,S_t(\omega),Y_t(\omega))$.

\begin{lemma} For all $K>0$, $ \Lambda_K$ is lower hemicontinuous, where 
$\Lambda_K (t,s,y)=\{(\alpha, \nu) \in \Lambda(t,s,y)|\; \alpha^*\alpha +\nu^*\nu \leq K\}$.
\label{lemmalhc}
\end{lemma}
{\bf Proof}
Let $(t_n,s_n,y_n)$ with limit $(t,s,y)$. Let $(\alpha,\nu) \in  \Lambda_K(t,s,y)$. Let $\nu_n=\nu \frac{\sqrt{(1-\rho(t_n,s_n,y_n)^2)}}{\sqrt{(1-\rho(t,s,y)^2)}}$ and $\alpha_n=-\nu\frac{\rho(t_n,s_n,y_n)}{\sqrt{(1-\rho(t,s,y)^2)}}$. 
It is easy to verify that $(\alpha_n,\nu_n)$ belongs to  $\Lambda(t_n,s_n,y_n)$, and that  $(\alpha_n,\nu_n) \rightarrow  (\alpha,\nu)$ as $n \rightarrow \infty$. Furthermore $\alpha_n^* \alpha_n+\nu_n^*\nu_n=\alpha^* \alpha+\nu^*\nu$. Thus $(\alpha_n,\nu_n) \in \Lambda_K(t_n,s_n,y_n)$. From the characterization of lower hemicontinuity given in Theorem 16.21 of \cite{AB}, this proves the result.
 \hfill $\square$\\
Given $\phi$ real valued denote  $\Lambda^{\phi}$ the multivalued Borel mapping
\begin{equation}
\Lambda^{\phi}(t,s,y)=\{(\alpha,\nu) \in \Lambda(t,s,y)\; |\;\alpha^*\alpha +\nu^*\nu \leq \phi(s,y)\}
\label{eqlambphi}
\end{equation}

\begin{lemma}
 Let $\phi(s,y)=C(1+(s^*s+y^*y)^{\frac{1}{2}}$ for $C>0$ and $k \in \N$. Then $\Lambda^{\phi}$ satisfies hypothesis $\overline H_{\Lambda}$
\label{corlambphi}
\end{lemma}
{\bf Proof}
It is enough to  prove that for all $K>0$, $\Lambda^{\phi}_K$ is lower hemicontinuous. 
 The proof is adapted from the proof of Lemma \ref{lemmalhc}, replacing $\nu_n$ and $\alpha_n$ respectively by $\tilde \nu_n=\nu_n \inf \big( 1,\sqrt{\frac{\phi(s_n, y_n)}{\phi(s, y)}} \big)$ and $\tilde \alpha_n= \alpha_n \inf \big( 1,\sqrt{\frac{\phi(s_n, y_n)}{\phi(s, y)}} \big)$.
 \hfill $\square$
The following proposition is an  application of our above results and of \cite{BN05} Theorem 5.1.
\begin{proposition} 
\begin{itemize}
\item  $$\Pi^{0,(S_0,Y_0)}_{s,t}(\xi)=\esssup_{P^{\Sigma,\mu}_{0,(S_0,Y_0)} \in (\tilde{\cal Q}_{0,(S_0,Y_0)})_S (\Lambda)} E_{P^{\Sigma,\mu}_{0,(S_0,Y_0)}}(\xi|{\cal B}_{s})$$ is the   surreplication price.
\item Choose for $g$ a non negative function on $\R_+ \times \R^{2n} \times \R^{2n}$ such that $g(t,(s,y),(0,0))=0$ for all $t$. 
Then 
 $$\Pi^{0,(S_0,Y_0)}_{s,t}(\xi)=\esssup_{P^{\Sigma,\mu}_{0,(S_0,Y_0)} \in (\tilde{\cal Q}_{0,(S_0,Y_0)})_S (\Lambda^{\phi})} (E_{P^{\Sigma,\mu}_{0,(S_0,Y_0)}}(\xi|{\cal B}_{s})- \alpha_{s,t}(P^{\Sigma,\mu}_{0,(S_0,Y_0)})$$  
defines a No Free Lunch Time Consistent  Convex Pricing Procedure extending the dynamics of the asset $S_t$. $\Pi^{0,(S_0,Y_0)}_{s,t}(\xi)$ represents the dynamic ask price for the asset $\xi$,  $-\Pi^{0,(S_0,Y_0)}_{s,t}(-\xi)$ represents its dynamic bid price. 
\item The dynamic ask price $\Pi^{0,(S_0,Y_0)}_{s,t}(h(S_t))$ of an european option, for $h$ Borelian bounded from below is equal to $\overline h(s,X_s,Y_s)$ for a lower semi continuous function  $\overline h$  on $[0,t[ \times \R^{2n}$.  When $h$ is continuous and when the hypothesis of Theorem \ref{thmuniqvis} are satisfied, $\overline h$ is a viscosity  solution of the second order PDE (\ref{eqedp02}).
\end{itemize}
\end{proposition}
 The choice $g\geq 0$ non equal to $0$ allows to take into account the liquidity risk.  As a simple example we can make the choice of a power penalty, i.e. $g(t,s,y,\alpha,\nu)=\frac{||(\alpha,\nu)||^p}{p}$ for some $1<p<\infty$, if $(\alpha,\nu) \in \Lambda^{\phi}(t,s,y)$, and $ g(t,s,y,\alpha,\nu)=+\infty$ otherwise.  Let $\psi_{t,s,y}(z,z')=||z-\frac{\rho(t,s,y)}{\sqrt{(1-\rho(t,s,y)^2)}}z'||$.
Let $q$ be the conjugate exponent of $p$.  The function $f$ is then given on $\R_+ \times \R^{2n} \times \R^{2n}$ by $f(t,s,y,z,z')=\frac{\psi_{t,s,y}(z,z')^q}{q}$  when $\psi_{t,s,y}(z,z')\leq \Phi(t,s,y)^{\frac{p-1}{2}}$, and $f(t,s,y,z,z')=\psi_{t,s,y}(z,z')\Phi(t,s,y)^{\frac{1}{2}} -\frac{\Phi(t,s,y)^{\frac{p}{2}}}{p}$ otherwise. It  is easy to verify that $g$ and $f$ satisfy the hypothesis of Theorems  \ref{thmsuper} and \ref{thmsub}. When $h$ is H\"older continuous, this provides a viscosity supersolution and a viscosity subsolution of (\ref{eqedp02}).\\
We can also choose the function $g$ such that the dynamic bid price  $-\Pi_{st}(-h_i(S_{t_i}))$  and ask price $\Pi_{st}(h_i(S_{t_i}))$ for a finite number  of  given options $h_i(S_{t_i})$ are compatible with limit order books observed in the financial market for these options (\cite{BN06}).

\section{Appendix: Time consistent dynamic procedure}
Motivated by the problematic of the evaluation of financial risks, the notion of risk measures has been introduced at the begining of the 21 st century.
Here we give the definition of a time consistent dynamic risk procedure,   which is up to a minus sign a time consistent dynamic risk measure. 
In the following ${\cal T}$ can be  either the set of deterministic times or the set of stopping times taking a finite number of finite values or the set of all stopping times.
\begin{definition}
A ${\cal T}$-time consistent dynamic procedure   on a filtered probability space $(\Omega,{\cal F}_{\infty},({\cal F}_t)_{t \in \R^+},P)$ is a family $(\Pi_{s, t})_{0 \leq s \leq t}$, $s,t$ in ${\cal T}$: \\ 
$\Pi_{s, t}: L^{\infty}({\cal F}_{t}) \rightarrow L^{\infty}({\cal F}_{s})$ satisfying  
\begin{enumerate}
\item  For every $s \leq t$ the following four properties: 
\begin{itemize}
\item  monotonicity: $\;\;if X\leq Y \;\; then \;\;
\Pi_{s,t}(X) \leq \Pi_{s,t}(Y)$  
\item  translation invariance: $\forall \; Z \in L^{\infty}({\cal F}_{s})$
$$\Pi_{s,t}(X+Z)=\Pi_{s,t}(X)+Z$$
\item  convexity: $\forall \lambda \in [0,1]$\
$$\Pi_{s,t}(\lambda X +(1-\lambda)Y) \leq \lambda \Pi_{s,t}(X)+(1-\lambda) \Pi_{s,t}(Y)$$
\item  continuity from below: 
 for every increasing 
sequence $X_n$ of elements of  $L^{\infty}({\cal F}_{t})$ such that $X=\lim\;X_n$, 
the increasing 
sequence $\Pi_{s,t}(X_n) $ has the limit $\Pi_{s,t}(X)$.
\end{itemize}
\item time consistency: For $0  \leq r \leq s \leq t$  in ${\cal T}$,  for all  $X \in L^{\infty}({\cal F}_{t})$ 
$$\Pi_{r,s}(\Pi_{s,t}(X))=\Pi_{r,t}(X)$$ 
\end{enumerate}
\label{TCP}
\end{definition}

\begin{definition} 
The  dynamic procedure can have additional properties.
\begin{itemize}
\item  It is normalized if  $\forall s \leq t\;\;\Pi_{s,t}(0)=0$
\item It  is  sublinear if 
$$\forall s \leq t \;\forall \lambda>0\;\; \forall X \in  L^{\infty}({\cal F}_{t})\;\; 
\Pi_{s,t}(\lambda X)= \lambda \Pi_{s,t}(X)$$
\end{itemize}
\label{def2}
\end{definition}

The most important way of constructing time consistent dynamic processes is to construct a stable set of equivalent probability measures and to define on this set  a penalty which is local and satisfies the cocycle condition. Recall the following definitions introduced in \cite{BN03}.
\begin{definition}
 A set ${\cal Q}$ of equivalent probability measures  is ${\cal T}$-stable if it satisfies the two following properties:
\begin{enumerate}
\item Stability by composition\\
 For all $s $ in ${\cal T}$,  for all $Q$ and $R$  in ${\cal Q}$, there is a probability measure $S$ in ${\cal Q}$ such that 
\begin{equation}
 (\frac{dS}{dP})= \frac{(\frac{dQ}{dP})}{(\frac{dQ}{dP})_{s}} (\frac{dR}{dP})_{s}\;\; 
\label{eqcomp}
\end{equation}
where  $(\frac{dR}{dP})_{s}$ means $E(\frac{dR}{dP}|{\cal F}_{s})$.
\item Stability by bifurcation\\
 For all $s $ in ${\cal T}$,  for all  $Q$ and $R$  in ${\cal Q}$, for all  $A \in {\cal F}_s$ there is a probability measure $S$ in ${\cal Q}$ such that  for all $X \in L^{\infty}(\Omega, {\cal F}_T,P)$, 
\begin{equation}
E_S(X|{\cal F}_s)=1_A E_Q(X|{\cal F}_s)+ 1_{A^c} E_R(X|{\cal F}_s)
\label{eqbif}
\end{equation}
\end{enumerate}
\label{defstable}
\end{definition}

Delbaen introduced in  \cite{D}  the notion of m-stability which is the condition of stability by composition for stopping times taking both finite and infinite values.  Notice that the m-stability implies from Proposition 1 of   \cite{D} the stability by bifurcation.  

  Recall the following definition of a  penalty \cite{BN04}, and of the cocycle condition

\begin{definition} 
 A penalty function $\alpha$ defined on a ${\cal T}$-stable set ${\cal Q}$ of probability measures all equivalent is a family of maps $(\alpha_{s,t}),\; s \leq t $ in ${\cal T}$,  defined on ${\cal Q}$ with values in the set of ${\cal F}_{s}$-measurable maps.
\begin {itemize}
\item[i)] It is local: \\
if for all $Q,R$ in ${\cal Q}$, for all  $s$ in ${\cal T}$, for all $A$ in ${\cal F}_{s}$, the assertion $1_A E_Q(X|{\cal F}_{s})= 1_A E_R(X|{\cal F}_{s})$ for all $X$ in 
$L^{\infty}({\cal F}_{t})$ implies that $1_A \alpha_{s,t}(Q)=1_A \alpha_{s,t}(R)$.
\item[ii)] It satisfies the cocycle condition if  
for  $r \leq s \leq t$ in ${\cal T}$,  for all $Q$ in ${\cal Q}$,
$$\alpha_{r,t}(Q)=\alpha_{r,s }(Q)+E_Q(\alpha_{s,t}(Q)|{\cal F}_{r})\; $$
\end{itemize}
\label{defpen}
\end{definition}

\end{document}